\documentclass[a4paper,11pt,leqno]{article}
\usepackage{amsfonts,amssymb,amsmath,amsthm,latexsym}
\usepackage{mathrsfs}
\usepackage{fullpage}
\usepackage{epic,eepic}

\newtheorem{theorem}{Theorem}[section]
\newtheorem{proposition}[theorem]{Proposition}
\newtheorem{lemma}[theorem]{Lemma}
\newtheorem{corollary}[theorem]{Corollary}
\theoremstyle{definition}
\newtheorem{remark}[theorem]{Remark}

\numberwithin{equation}{section}

\newcommand{\N}{{\mathbb N}}
\newcommand{\R}{{\mathbb R}}
\newcommand{\E}{{\mathcal E}}
\newcommand{\eps}{{\varepsilon}}
\newcommand{\spt}{\mathrm{supp}}

\begin{document}

\title{Positive solutions to nonlinear $p$-Laplace equations\\
with Hardy potential in exterior domains}

\author{Vitali Liskevich, Sofya Lyakhova and Vitaly Moroz\medskip\\
School of Mathematics\\
University of Bristol\\
Bristol BS8 1TW\\
United Kingdom\smallskip\\
{\small E-mail: {\tt V.Liskevich$\,|\,$S.Lyakhova$\,|\,$V.Moroz@bristol.ac.uk}}}
\date{\today}

\maketitle

\begin{abstract}
We study the existence and nonexistence of positive
(super)\,solutions to the nonlinear $p$-Laplace equation
$$-\Delta_p u-\frac{\mu}{|x|^p}u^{p-1}=\frac{C}{|x|^{\sigma}}u^q$$
in exterior domains of ${\R}^N$ ($N\ge 2$). Here $p\in(1,+\infty)$
and $\mu\le C_H$, where $C_H$ is the critical Hardy constant.
We provide a sharp characterization of the set of $(q,\sigma)\in\R^2$
such that the equation has no positive (super) solutions.

The proofs are based on the explicit construction of appropriate
barriers and involve the analysis of asymptotic behavior of
super-harmonic functions associated to the $p$-Laplace operator
with Hardy-type potentials, comparison principles and an improved
version of Hardy's inequality in exterior domains.  In the context of the $p$-Laplacian
we establish the existence and asymptotic behavior of the
harmonic functions by means of the generalized Pr\"ufer-Transformation.
\end{abstract}

\begin{small}
\tableofcontents
\end{small}

\section{Introduction and Results}
We study the problem of the existence and nonexistence of positive
(super)\,solutions to nonlinear $p$-Laplace equation with Hardy potential
\begin{equation}\label{e:MAIN}
-\Delta_p u
-\frac{\mu}{|x|^p}u^{p-1}=\frac{C}{|x|^\sigma}u^q\quad\mbox{in}\quad B^c_\rho,
\end{equation}
where $-\Delta_p u=-\mathrm{div}(|\nabla u|^{p-2}\nabla u)$ is the $p$-Laplace operator,
$1<p<\infty$, $C>0$, $\mu\in\R$, $(q,\sigma)\in\R^2$ and
$B^c_\rho:=\{x\in\R^N:\,|x|>\rho\}$ is the exterior of the ball in $\R^N$, with $N\geq 2$.
We say that $u\in W^{1,p}_{loc}(G)\cap C(G)$ is a  \textit{super-solution} to equation \eqref{e:MAIN} in a domain
$G\subseteq\R^N$ with $0\not\in G$ if for all
$0\leq\varphi\in W^{1,p}_c(G)\cap C(G)$ the following inequality holds
$$\int_G\nabla u|\nabla u|^{p-2}\nabla\varphi\,dx-
\int_G\frac{\mu}{|x|^p}u^{p-1} \varphi\, dx
\geq\int_G\frac{C}{|x|^\sigma}u^{q}\varphi\, dx.$$  Here and below
$W^{1,p}_c(G):=\{u\in W^{1,p}_{loc}(G),\,\spt(u)\Subset G\}$.
The notions of a sub-solution and solution are defined similarly,
by replacing $"\geq"$ with $"\leq"$ and $"="$, respectively.
It follows from the Harnack inequality  (cf. \cite{Serrin64})
that any nontrivial nonnegative super-solution to \eqref{e:MAIN} in $G$
is strictly positive in $G$.

One of the features of equation \eqref{e:MAIN} on unbounded
domains is the nonexistence of positive solutions for certain
values of the exponent $q$. Such Liouville type nonexistence phenomena
have been known for semilinear elliptic equations ($p=2$) at least since
the celebrated  works of Serrin in the earlier 70's (cf. the references in \cite{Serrin-Zou})
and of Gidas and Spruck \cite{Gidas-Spruck}.
One of the first Liouville-type results for the nonlinear $p$-Laplace equations in exterior domains
is due to Bidaut--V\'eron \cite[Theorem 1.3]{B-Veron}. Theorem A below extends the result in \cite{B-Veron},
including the cases $p>N$ and $q<p-1$.

\smallskip\noindent
\textbf{Theorem A\,} (\cite[Theorem 1.3]{B-Veron}, Theorem \ref{t:MAIN} below).
\textit{The equation
\begin{equation}\label{e-BV}
-\Delta_p u=u^q\quad\mbox{in}\quad B^c_\rho
\end{equation}
has no positive super-solutions if and only if  $q_\ast\le q \le q^\ast$, where $q^\ast=\frac{N(p-1)}{N-p}$ when $p<N$,
or $q^\ast=+\infty$ when $p\ge N$, and $q_\ast=-\infty$ when $p\le N$ or $q^\ast=\frac{N(p-1)}{N-p}$ when $p>N$.}
\smallskip

Theorem A had been generalized and extended in various direction by many authors
(see, e.g., \cite{Peral-1,B-VP,Pohozhaev,Serrin-Zou,Veron} and  references therein).
The techniques in those works usually involve careful integral estimates
and/or sophisticated analysis of related nonlinear ODE's.
A different approach to nonlinear Liouville type theorems goes to back to an earlier paper by Kondratiev and Landis
\cite{Kondratiev-Landis} and was recently developed in the context of semilinear equations ($p=2$)
in \cite{KLS,KLS1,KLM,KLMS}.
The approach is based on the pointwise Phragm\'en--Lindel\"of type bounds
on positive super-harmonic functions and related Hardy--type inequalities.

Recall that the classical {\em Hardy inequality} states that
\begin{equation}\label{e:HARDY-0}
\int_{B_\rho^c}|\nabla u|^pdx\geq
C_H\int_{B_\rho^c}\frac{u^p}{|x|^{p}}dx,\qquad\forall\:u\in C^\infty_c(B_\rho^c),
\end{equation}
with the sharp constant $C_H=\left|\frac{N-p}{p}\right|^p$, $p>1$.
The optimality of the Hardy constant $C_H$ implies, via Picone's identity, the following nonexistence result.

\smallskip\noindent
\textbf{Theorem B\,} (\cite[Corollary 2.2]{AllegrettoHuang}, Corollary \ref{c:NonExist} below)
\textit{The equation
\begin{eqnarray}\label{e:MAIN-0}
-\Delta_p u-\frac{\mu}{|x|^p}u^{p-1}=0\quad\mbox{in}\quad B^c_\rho
\end{eqnarray}
has no positive super-solutions if and only if $\mu>C_H$.}
\smallskip

Let us sketch a simple proof of the nonexistence part of Theorem A in the case $p\neq N$ and $q_\ast<q<q^\ast$.
Indeed, let $u>0$ be a  super-solution to \eqref{e-BV}.
Then $-\Delta_p u\ge 0$ in $B_\rho^c$.
A comparison principle for the $p$-Laplacian in exterior domains (see Theorem \ref{l:WCPimpr}
and Theorems \ref{l:WCP} and \ref{l:upper_est<C_H} below) implies that $u$  obeys
the {\em Phragm\'en--Lindel\"of type bounds}
\begin{equation}\label{BV-1}
c|x|^{\gamma_-}\le\inf_{|x|=r}u\le c^{-1}|x|^{\gamma_+}\quad\mbox{in}\quad B^c_{2\rho},
\end{equation}
where $\gamma_-=\min\{0,\frac{p-N}{p-1}\}$ and $\gamma_+=\max\{0,\frac{p-N}{p-1}\}$.

Assume that $q\ge p-1$ and perform a {\em homogenization} of equation \eqref{e-BV} rewriting it in the from
\begin{equation}\label{BV-hom}
-\Delta_p u=V(x)u^{p-1}\quad\mbox{in}\quad B^c_\rho,
\end{equation}
where $V(x):=u^{q-(p-1)}$. Using the lower bound from \eqref{BV-1}, we conclude that
$$V(x)\ge c_1|x|^{\gamma_-(q-p+1)}\quad\mbox{in }\:B_{2\rho}^c.$$
Hence, by Theorem B, \eqref{BV-hom} has no positive super-solutions provided that $\gamma_-(q-p+1)>-p$.
Therefore \eqref{e-BV} has no positive super-solutions when $p-1\le q < q^\ast$.

Now assume that $q<p-1$. Then a standard scaling argument (see  Lemma \ref{l:EstSubLin} below)
shows that any super-solution $u>0$ to nonlinear equation \eqref{e-BV} obeys the lower bound
$$u\ge c|x|^{\frac{p}{(p-1)-q}}\quad\text{in }\: \:B_{2\rho}^c.$$
Comparing this "nonlinear" estimate with the upper bound in \eqref{BV-1}, we conclude that
equation \eqref{e-BV} has no positive super-solutions for $q_\ast<q<p-1$.

The above simple proof relies only on Theorem B and pointwise Phragm\'en--Lindel\"of type bounds \eqref{BV-1}.
It does not cover the {\sl critical cases} $q=q_\ast$ and $q=q^\ast$,
where additional arguments are required. On the other hand,
an explicit construction of radial super-solutions to \eqref{e-BV}
when $q\not\in[q_\ast,q^\ast]$ shows that the values of the {\sl critical exponents}
$q^\ast$ and $q_\ast$ are sharp.
Considerations of this type first appeared in \cite{KLS}. They have proved to be a
powerful and flexible tool for studying nonlinear Liouville phenomena for various classes of elliptic
operators and domains, see \cite{KLS,KLS1,KLM,KLMS,LLM,LiskevichSkrypnikSkrypnik}.

In this paper we are interested in nonlinear Liouville theorems for perturbations
of the $p$-Laplace operator by the Hardy type potential.
To explore the impact of the potential on the value of the critical exponents $q^\ast$ and $q_\ast$, let us
consider the equation of the form
\begin{equation}\label{e:main_e}
-\Delta_p u-\frac{\mu}{|x|^{p+\epsilon}}u^{p-1}=u^q\quad \mbox{in }\:
B^c_\rho,
\end{equation}
where $\mu\in\R$ and $\epsilon\in\R$.
One can verify directly that
{\sl if $\epsilon<0$ and $\mu<0$, then \eqref{e:main_e} admits positive solutions for all $q\in\R$},
while
{\sl
if $\epsilon<0$ and $\mu>0$ then \eqref{e:main_e} has no positive super-solutions for any
$q\in\R$}. The latter follows immediately from Theorem B.
On the other hand, one can show (see \cite[Theorem 1.2]{LiskevichSkrypnikSkrypnik}) that if $\epsilon>0$ then
\eqref{e:main_e} {\sl has the same critical exponents $q^\ast$ and $q_\ast$ as \eqref{e-BV}}.
This follows from the fact that positive super-solutions to
\begin{equation}\label{e-V-eps}
-\Delta_p u-\frac{\mu}{|x|^{p+\epsilon}}u^{p-1}=0\quad \mbox{in }\: B^c_\rho
\end{equation}
satisfy the same bound \eqref{BV-1} as super-solutions to $-\Delta_p u\ge 0$ in $B^c_\rho$.

In this paper we show that in the borderline case $\epsilon=0$ the
critical exponents for equation \eqref{e:main_e} explicitly depend on $\mu$.
This is a consequence of the fact that the Phragm\'en--Lindel\"of bounds
for equation \eqref{e-V-eps} with $\epsilon=0$
become sensitive to the value of the parameter $\mu$.
Such phenomenon and its relation to the Hardy
type inequalities has been recently observed  for $p=2$ in the
case of the ball as well as exterior domains in
\cite{Brezis-Cabre,Brezis-Tesei,Dupaigne,LLM,Pohozaev-Tesei,Smets,Terracini}. The main
difficulty comparing with the semilinear case $p=2$ arises when a
comparison principle for the $p$--Laplacian has to be involved in
the argument. After examples in \cite{DelPino,Fleck-1} it is known that solutions to the equation
$-\Delta_p u-V(x)u^{p-1}=0$ may not satisfy the usual comparison
principle as soon as the potential has a nontrivial negative part
$(V^+)$. The proof of the (restricted) comparison principle
requires delicate arguments. We provide a new version of the
comparison principle (Theorem \ref{l:WCPimpr}), following  the
ideas from \cite{ShafrirMarcus}. In order to use this result for
obtaining sharp Phragm\'en--Lindel\"of bounds one has to produce
explicitly a radial sub-solution to a homogeneous equation in the
exterior of the ball with zero data on the sphere. This has been
resolved in this paper by means of the generalized Pr\"ufer
transformation (see \cite{ReichelWalter} and Appendix B.3), which,
up to our knowledge, has never been used before in this context.
We also provide an elementary proof of an improved Hardy
Inequality in exterior domains.
Improved Hardy Inequality plays a crucial role in our
analysis of equation \eqref{e:MAIN} in the critical case
$\mu=C_H$.

To formulate the main result of the paper we assume that $\mu\leq C_H$,
otherwise \eqref{e:MAIN} has no positive super-solutions by Theorem B.
When $\mu\leq C_H$, the scalar equation
\begin{eqnarray*}
-\gamma|\gamma|^{p-2}\left(\gamma(p-1)+N-p\right)=\mu
\end{eqnarray*}
has two real roots $\gamma_-\leq\gamma_+$. Note that
if $\mu=C_H$ then $\gamma_-=\gamma_+=\frac{p-N}{p}$. For $\mu\leq
C_H$ we introduce the critical line $\Lambda^*(q,\mu)$ for
equation \eqref{e:MAIN} on the $(q,\sigma)$--plane
$$\Lambda_*(q,\mu):=\min\{\gamma_-(q-p+1)+p,\;\gamma_+(q-p+1)+p\}\qquad (q\in\R),$$
and the nonexistence set
$$\mathcal{N}=\{(q,\sigma)\in\R^2\setminus(p-1,p):
\textit{\eqref{e:MAIN} has no positive supersolutions  in } B^c_\rho\}.$$

\begin{theorem}\label{t:MAIN}
The following assertions are valid.
\begin{enumerate}
\item[{$(i)$}]
If $\mu<C_H$ then $\mathcal N=\{\sigma\leq \Lambda_*(q)\}$.
\item[{$(ii)$}]
If $\mu=C_H$ then $\mathcal
N=\{\sigma<\Lambda_\ast(q)\}\cup\{\sigma=\Lambda_\ast(q),\,q\geq -1\}$.
\end{enumerate}
\end{theorem}

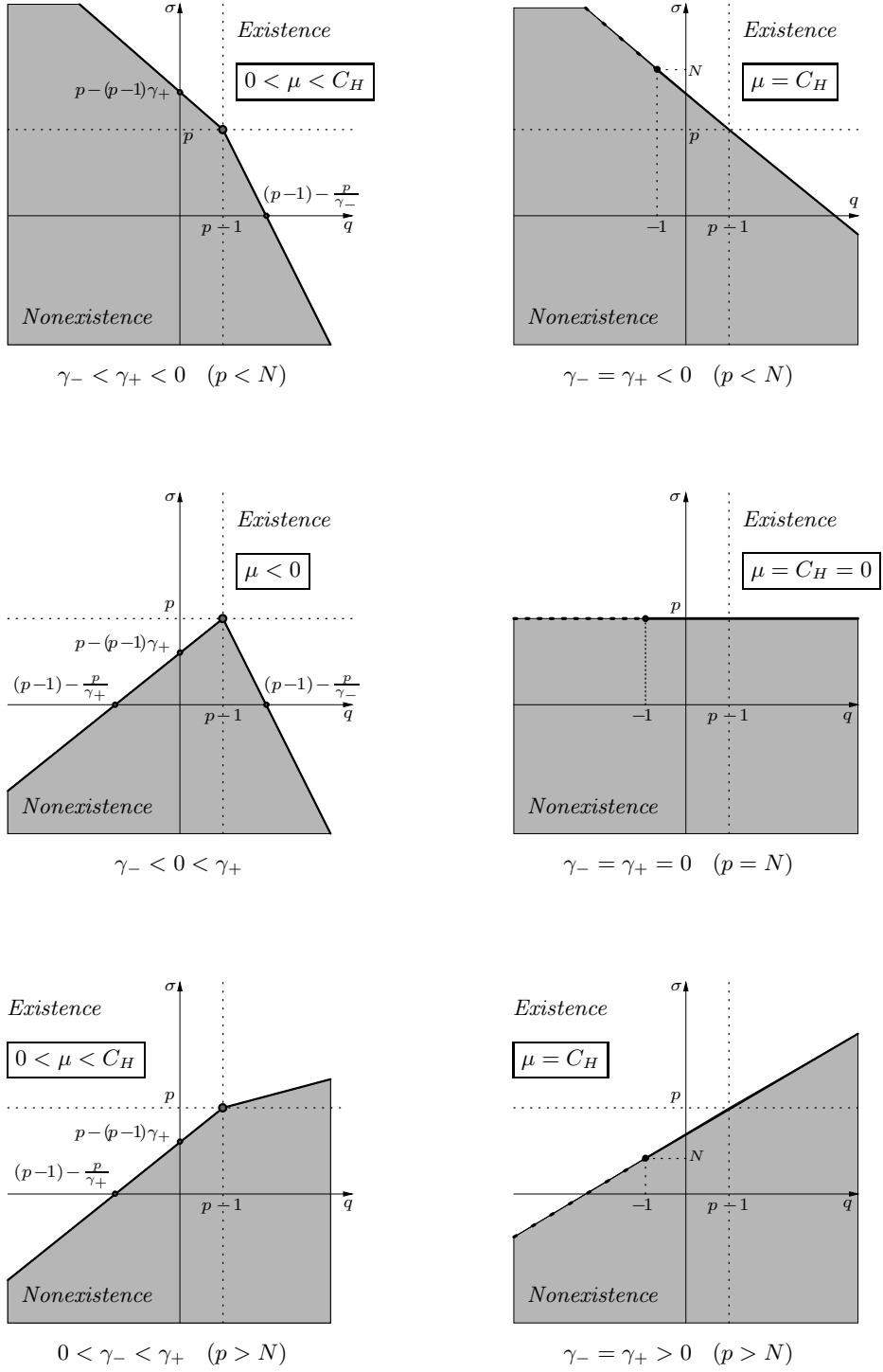
\begin{figure}[p]
\begin{center}

\setlength{\unitlength}{0.0004in}
\begingroup\makeatletter\ifx\SetFigFont\undefined%
\gdef\SetFigFont#1#2#3#4#5{%
  \reset@font\fontsize{#1}{#2pt}%
  \fontfamily{#3}\fontseries{#4}\fontshape{#5}%
  \selectfont}%
\fi\endgroup%
{\renewcommand{\dashlinestretch}{30}


\vspace{0.3cm}


\begin{picture}(4815,4800)(0,-300) 
\put(600,-500){\makebox(0,0)[lb]{\smash{{{\SetFigFont{9}{12.0}{\rmdefault}{\mddefault}{\updefault}
$\gamma_-<\gamma_+<0\;\;$ ($p<N$)}}}}} \texture{88555555 55000000
555555 55000000 555555 55000000 555555 55000000
           555555 55000000 555555 55000000 555555 55000000 555555 55000000
           555555 55000000 555555 55000000 555555 55000000 555555 55000000
           555555 55000000 555555 55000000 555555 55000000 555555 55000000 }
\shade\path(3012,3012)(4512,12)(12,12)(12,4762)(1012,4762)(3012,3012)
\shade\path(3012,3012)(4512,12)(12,12)(12,4762)(1012,4762)(3012,3012)

\path(12,1812)(4812,1812) 
\blacken\path(4692.000,1782.000)(4812.000,1812.000)(4692.000,1842.000)(4692.000,1782.000) 
\path(2412,12)(2412,4762) 
\blacken\path(2442.000,4642.000)(2412.000,4762.000)(2382.000,4642.000)(2442.000,4642.000) 

\dottedline{120}(12,3012)(4812,3012)
\dottedline{120}(3012,4762)(3012,12)

\thicklines \path(1012,4762)(3012,3012)(4512,12)

\texture{aa777777 77aaaaaa aad5d5d5 d5aaaaaa aa777777 77aaaaaa
aadd5ddd 5daaaaaa
         aa777777 77aaaaaa aad5d5d5 d5aaaaaa aa777777 77aaaaaa aa5ddd5d ddaaaaaa
         aa777777 77aaaaaa aad5d5d5 d5aaaaaa aa777777 77aaaaaa aadd5ddd 5daaaaaa
         aa777777 77aaaaaa aad5d5d5 d5aaaaaa aa777777 77aaaaaa aa5ddd5d ddaaaaaa }
\put(3012,3012){\shade\ellipse{100}{100}}
\put(3620,1812){\shade\ellipse{60}{60}}
\put(2412,3532){\shade\ellipse{60}{60}}

\put(4690,1630){\makebox(0,0)[lb]{\smash{{{\SetFigFont{7}{12.0}{\rmdefault}{\mddefault}{\updefault}$q$}}}}}
\put(2200,4650){\makebox(0,0)[lb]{\smash{{{\SetFigFont{7}{12.0}{\rmdefault}{\mddefault}{\updefault}$\sigma$}}}}}
\put(2460,2870){\makebox(0,0)[lb]{\smash{{{\SetFigFont{7}{12.0}{\rmdefault}{\mddefault}{\updefault}$p$}}}}}
\put(2760,1590){\makebox(0,0)[lb]{\smash{{{\SetFigFont{7}{12.0}{\rmdefault}{\mddefault}{\updefault}$\!p-\!1$}}}}}
\put(3570,2050){\makebox(0,0)[lb]{\smash{{{\SetFigFont{7}{12.0}{\rmdefault}{\mddefault}{\updefault}$(p\!-\!\!1)\!-\!\frac{p}{\gamma_-}$}}}}}
\put(950,3500){\makebox(0,0)[lb]{\smash{{{\SetFigFont{7}{12.0}{\rmdefault}{\mddefault}{\updefault}$p\!-\!(\!p\!-\!\!1\!)\gamma_+$}}}}}
\put(200,300){\makebox(0,0)[lb]{\smash{{{\SetFigFont{9}{12.0}{\rmdefault}{\mddefault}{\itdefault}Nonexistence}}}}}
\put(3200,4300){\makebox(0,0)[lb]{\smash{{{\SetFigFont{9}{12.0}{\rmdefault}{\mddefault}{\itdefault}Existence}}}}}
\put(3200,3600){\makebox(0,0)[lb]{\smash{{{\SetFigFont{9}{12.0}{\rmdefault}{\mddefault}{\updefault}\boxed{0<\mu<C_H}}}}}}
\end{picture}
\hspace{2cm}
\begin{picture}(4815,4800)(0,-300) 
\put(600,-500){\makebox(0,0)[lb]{\smash{{{\SetFigFont{9}{12.0}{\rmdefault}{\mddefault}{\updefault}
$\gamma_-=\gamma_+<0\;\;$ ($p<N$)}}}}}
\texture{88555555 55000000
555555 55000000 555555 55000000 555555 55000000
           555555 55000000 555555 55000000 555555 55000000 555555 55000000
           555555 55000000 555555 55000000 555555 55000000 555555 55000000
           555555 55000000 555555 55000000 555555 55000000 555555 55000000 }
\shade\path(3012,3012)(4812,1552)(4812,12)(12,12)(12,4712)(1012,4712)(3012,3012)

\path(12,1812)(4812,1812) 
\blacken\path(4692.000,1782.000)(4812.000,1812.000)(4692.000,1842.000)(4692.000,1782.000) 
\path(2412,12)(2412,4762) 
\blacken\path(2442.000,4642.000)(2412.000,4762.000)(2382.000,4642.000)(2442.000,4642.000) 

\dottedline{120}(12,3012)(4812,3012)
\dottedline{120}(3012,4762)(3012,12)
\dottedline{120}(2012,3852)(2012,1812)
\dottedline{120}(2012,3852)(2412,3852)

\thicklines
\path(2012,3852)(3012,3012)(4812,1552)



\texture{aa777777 77aaaaaa aad5d5d5 d5aaaaaa aa777777 77aaaaaa
aadd5ddd 5daaaaaa
         aa777777 77aaaaaa aad5d5d5 d5aaaaaa aa777777 77aaaaaa aa5ddd5d ddaaaaaa
         aa777777 77aaaaaa aad5d5d5 d5aaaaaa aa777777 77aaaaaa aadd5ddd 5daaaaaa
         aa777777 77aaaaaa aad5d5d5 d5aaaaaa aa777777 77aaaaaa aa5ddd5d ddaaaaaa }
\Thicklines
\dottedline{300}(1012,4712)(2012,3852)
\put(2012,3852){\ellipse{60}{60}}

\put(4690,1960){\makebox(0,0)[lb]{\smash{{{\SetFigFont{7}{12.0}{\rmdefault}{\mddefault}{\updefault}$q$}}}}}
\put(2200,4650){\makebox(0,0)[lb]{\smash{{{\SetFigFont{7}{12.0}{\rmdefault}{\mddefault}{\updefault}$\sigma$}}}}}
\put(2460,2870){\makebox(0,0)[lb]{\smash{{{\SetFigFont{7}{12.0}{\rmdefault}{\mddefault}{\updefault}$p$}}}}}
\put(2480,3770){\makebox(0,0)[lb]{\smash{{{\SetFigFont{6}{12.0}{\rmdefault}{\mddefault}{\updefault}$\!N$}}}}}

\put(2760,1590){\makebox(0,0)[lb]{\smash{{{\SetFigFont{7}{12.0}{\rmdefault}{\mddefault}{\updefault}$\!p-\!1$}}}}}
\put(1850,1590){\makebox(0,0)[lb]{\smash{{{\SetFigFont{7}{12.0}{\rmdefault}{\mddefault}{\updefault}$-\!1$}}}}}


\put(200,300){\makebox(0,0)[lb]{\smash{{{\SetFigFont{9}{12.0}{\rmdefault}{\mddefault}{\itdefault}Nonexistence}}}}}
\put(3200,4300){\makebox(0,0)[lb]{\smash{{{\SetFigFont{9}{12.0}{\rmdefault}{\mddefault}{\itdefault}Existence}}}}}

\put(3200,3600){\makebox(0,0)[lb]{\smash{{{\SetFigFont{9}{12.0}{\rmdefault}{\mddefault}{\itdefault}$\boxed{\mu=C_H}$}}}}}

\end{picture}
\vspace{2cm}


\begin{picture}(4815,4800)(0,-300) 
\texture{88555555 55000000 555555 55000000 555555 55000000 555555
55000000
           555555 55000000 555555 55000000 555555 55000000 555555 55000000
           555555 55000000 555555 55000000 555555 55000000 555555 55000000
           555555 55000000 555555 55000000 555555 55000000 555555 55000000 }
\shade\path(3012,3012)(4512,12)(12,12)(12,612)(3012,3012)
\shade\path(3012,3012)(4512,12)(12,12)(12,612)(3012,3012)

\path(12,1812)(4812,1812) 
\blacken\path(4692.000,1782.000)(4812.000,1812.000)(4692.000,1842.000)(4692.000,1782.000) 
\path(2412,12)(2412,4642) 
\blacken\path(2442.000,4642.000)(2412.000,4762.000)(2382.000,4642.000)(2442.000,4642.000) 

\dottedline{120}(12,3012)(4812,3012)
\dottedline{120}(3012,4762)(3012,12) \thicklines
\path(12,612)(3012,3012)(4512,12) 

\texture{aa777777 77aaaaaa aad5d5d5 d5aaaaaa aa777777 77aaaaaa
aadd5ddd 5daaaaaa
         aa777777 77aaaaaa aad5d5d5 d5aaaaaa aa777777 77aaaaaa aa5ddd5d ddaaaaaa
         aa777777 77aaaaaa aad5d5d5 d5aaaaaa aa777777 77aaaaaa aadd5ddd 5daaaaaa
         aa777777 77aaaaaa aad5d5d5 d5aaaaaa aa777777 77aaaaaa aa5ddd5d ddaaaaaa }
\put(3012,3012){\shade\ellipse{100}{100}}
\put(1510,1812){\shade\ellipse{60}{60}}
\put(3620,1812){\shade\ellipse{60}{60}}
\put(2412,2540){\shade\ellipse{60}{60}}

\put(4690,1630){\makebox(0,0)[lb]{\smash{{{\SetFigFont{7}{12.0}{\rmdefault}{\mddefault}{\updefault}$q$}}}}}
\put(2200,4650){\makebox(0,0)[lb]{\smash{{{\SetFigFont{7}{12.0}{\rmdefault}{\mddefault}{\updefault}$\sigma$}}}}}
\put(2200,3170){\makebox(0,0)[lb]{\smash{{{\SetFigFont{7}{12.0}{\rmdefault}{\mddefault}{\updefault}$p$}}}}}
\put(2760,1590){\makebox(0,0)[lb]{\smash{{{\SetFigFont{7}{12.0}{\rmdefault}{\mddefault}{\updefault}$p\!-\!1$}}}}}
\put(3580,2030){\makebox(0,0)[lb]{\smash{{{\SetFigFont{7}{12.0}{\rmdefault}{\mddefault}{\updefault}$(p\!-\!\!1)\!-\!\frac{p}{\gamma_-}$}}}}}
\put(80,2030){\makebox(0,0)[lb]{\smash{{{\SetFigFont{7}{12.0}{\rmdefault}{\mddefault}{\updefault}$(p\!-\!\!1)\!-\!\frac{p}{\gamma_+}$}}}}}
\put(950,2600){\makebox(0,0)[lb]{\smash{{{\SetFigFont{7}{12.0}{\rmdefault}{\mddefault}{\updefault}$p\!-\!(\!p\!-\!\!1\!)\gamma_+$}}}}}
\put(200,300){\makebox(0,0)[lb]{\smash{{{\SetFigFont{9}{12.0}{\rmdefault}{\mddefault}{\itdefault}Nonexistence}}}}}
\put(3200,4300){\makebox(0,0)[lb]{\smash{{{\SetFigFont{9}{12.0}{\rmdefault}{\mddefault}{\itdefault}Existence}}}}}
\put(3200,3600){\makebox(0,0)[lb]{\smash{{{\SetFigFont{9}{12.0}{\rmdefault}{\mddefault}{\itdefault}$\boxed{\mu<0}$}}}}}
\put(1400,-500){\makebox(0,0)[lb]{\smash{{{\SetFigFont{9}{12.0}{\rmdefault}{\mddefault}{\updefault}
$\gamma_-<0<\gamma_+$}}}}}
\end{picture}
\hspace{2cm}
\begin{picture}(4815,4800)(0,-300) 

\texture{88555555 55000000 555555 55000000 555555 55000000 555555
55000000
           555555 55000000 555555 55000000 555555 55000000 555555 55000000
           555555 55000000 555555 55000000 555555 55000000 555555 55000000
           555555 55000000 555555 55000000 555555 55000000 555555 55000000 }
\thinlines\shade\path(3012,3012)(4812,3012)(4812,12)(12,12)(12,3012)(3012,3012)

\path(12,1812)(4812,1812) 
\blacken\path(4692.000,1782.000)(4812.000,1812.000)(4692.000,1842.000)(4692.000,1782.000) 
\path(2412,12)(2412,4762) 
\blacken\path(2442.000,4642.000)(2412.000,4762.000)(2382.000,4642.000)(2442.000,4642.000) 

\dottedline{120}(3012,4762)(3012,12)
\dottedline{60}(1850,3012)(1850,1812)

\Thicklines
\drawline(1850,3012)(4812,3012) 
\Thicklines
\dottedline{180}(12,3012)(1750,3012) 
\put(1850,3012){\ellipse{60}{60}}

\put(4600,1600){\makebox(0,0)[lb]{\smash{{{\SetFigFont{7}{12.0}{\rmdefault}{\mddefault}{\updefault}$q$}}}}}
\put(2180,4650){\makebox(0,0)[lb]{\smash{{{\SetFigFont{7}{12.0}{\rmdefault}{\mddefault}{\updefault}$\sigma$}}}}}
\put(2200,3130){\makebox(0,0)[lb]{\smash{{{\SetFigFont{7}{12.0}{\rmdefault}{\mddefault}{\updefault}$p$}}}}}
\put(2760,1590){\makebox(0,0)[lb]{\smash{{{\SetFigFont{7}{12.0}{\rmdefault}{\mddefault}{\updefault}$\!p-\!1$}}}}}
\put(1650,1590){\makebox(0,0)[lb]{\smash{{{\SetFigFont{7}{12.0}{\rmdefault}{\mddefault}{\updefault}$-\!1$}}}}}
\put(200,300){\makebox(0,0)[lb]{\smash{{{\SetFigFont{9}{12.0}{\rmdefault}{\mddefault}{\itdefault}Nonexistence}}}}}
\put(3200,4300){\makebox(0,0)[lb]{\smash{{{\SetFigFont{9}{12.0}{\rmdefault}{\mddefault}{\itdefault}Existence}}}}}
\put(3200,3600){\makebox(0,0)[lb]{\smash{{{\SetFigFont{9}{12.0}{\rmdefault}{\mddefault}{\itdefault}$\boxed{\mu=C_H=0}$}}}}}
\put(600,-500){\makebox(0,0)[lb]{\smash{{{\SetFigFont{9}{12.0}{\rmdefault}{\mddefault}{\updefault}
$\gamma_-=\gamma_+=0\;\;$ ($p=N$)}}}}}
\end{picture}
\vspace{2cm}


\begin{picture}(4815,4800)(0,-300) 
\put(600,-500){\makebox(0,0)[lb]{\smash{{{\SetFigFont{9}{12.0}{\rmdefault}{\mddefault}{\updefault}
$0<\gamma_-<\gamma_+\;\;$ ($p>N$)}}}}} \texture{88555555 55000000
555555 55000000 555555 55000000 555555 55000000
           555555 55000000 555555 55000000 555555 55000000 555555 55000000
           555555 55000000 555555 55000000 555555 55000000 555555 55000000
           555555 55000000 555555 55000000 555555 55000000 555555 55000000 }
\shade\path(3012,3012)(4512,3412)(4512,12)(12,12)(12,612)(3012,3012)
\shade\path(3012,3012)(4512,3412)(4512,12)(12,12)(12,612)(3012,3012)


\path(12,1812)(4812,1812) 
\blacken\path(4692.000,1782.000)(4812.000,1812.000)(4692.000,1842.000)(4692.000,1782.000) 
\path(2412,12)(2412,4762) 
\blacken\path(2442.000,4642.000)(2412.000,4762.000)(2382.000,4642.000)(2442.000,4642.000) 

\dottedline{120}(12,3012)(4650,3012)
\dottedline{120}(3012,4762)(3012,12)

\thicklines \path(12,612)(3012,3012)(4512,3412)

\texture{aa777777 77aaaaaa aad5d5d5 d5aaaaaa aa777777 77aaaaaa
aadd5ddd 5daaaaaa
         aa777777 77aaaaaa aad5d5d5 d5aaaaaa aa777777 77aaaaaa aa5ddd5d ddaaaaaa
         aa777777 77aaaaaa aad5d5d5 d5aaaaaa aa777777 77aaaaaa aadd5ddd 5daaaaaa
         aa777777 77aaaaaa aad5d5d5 d5aaaaaa aa777777 77aaaaaa aa5ddd5d ddaaaaaa }
\put(3012,3012){\shade\ellipse{100}{100}}
\put(1510,1812){\shade\ellipse{60}{60}}
\put(2412,2540){\shade\ellipse{60}{60}}

\put(4690,1630){\makebox(0,0)[lb]{\smash{{{\SetFigFont{7}{12.0}{\rmdefault}{\mddefault}{\updefault}$q$}}}}}
\put(2200,4650){\makebox(0,0)[lb]{\smash{{{\SetFigFont{7}{12.0}{\rmdefault}{\mddefault}{\updefault}$\sigma$}}}}}
\put(2200,3130){\makebox(0,0)[lb]{\smash{{{\SetFigFont{7}{12.0}{\rmdefault}{\mddefault}{\updefault}$p$}}}}}
\put(2760,1590){\makebox(0,0)[lb]{\smash{{{\SetFigFont{7}{12.0}{\rmdefault}{\mddefault}{\updefault}$\!p-\!1$}}}}}
\put(100,2050){\makebox(0,0)[lb]{\smash{{{\SetFigFont{7}{12.0}{\rmdefault}{\mddefault}{\updefault}$(p\!-\!\!1)\!-\!\frac{p}{\gamma_+}$}}}}}
\put(950,2600){\makebox(0,0)[lb]{\smash{{{\SetFigFont{7}{12.0}{\rmdefault}{\mddefault}{\updefault}$p\!-\!(\!p\!-\!\!1\!)\gamma_+$}}}}}
\put(200,300){\makebox(0,0)[lb]{\smash{{{\SetFigFont{9}{12.0}{\rmdefault}{\mddefault}{\itdefault}Nonexistence}}}}}
\put(0,4300){\makebox(0,0)[lb]{\smash{{{\SetFigFont{9}{12.0}{\rmdefault}{\mddefault}{\itdefault}Existence}}}}}
\put(0,3600){\makebox(0,0)[lb]{\smash{{{\SetFigFont{9}{12.0}{\rmdefault}{\mddefault}{\updefault}\boxed{0<\mu<C_H}}}}}}
\end{picture}
\hspace{2cm}
\begin{picture}(4815,4800)(0,-300) 
\put(600,-500){\makebox(0,0)[lb]{\smash{{{\SetFigFont{9}{12.0}{\rmdefault}{\mddefault}{\updefault}
$\gamma_-=\gamma_+>0\;\;$ ($p>N$)}}}}} \texture{88555555 55000000
555555 55000000 555555 55000000 555555 55000000
           555555 55000000 555555 55000000 555555 55000000 555555 55000000
           555555 55000000 555555 55000000 555555 55000000 555555 55000000
           555555 55000000 555555 55000000 555555 55000000 555555 55000000 }

\shade\path(3012,3012)(4812,4046)(4812,12)(12,12)(12,1212)(3012,3012)

\path(12,1812)(4812,1812) 
\blacken\path(4692.000,1782.000)(4812.000,1812.000)(4692.000,1842.000)(4692.000,1782.000) 
\path(2412,12)(2412,4762) 
\blacken\path(2442.000,4642.000)(2412.000,4762.000)(2382.000,4642.000)(2442.000,4642.000) 

\dottedline{120}(12,3012)(4812,3012)
\dottedline{120}(3012,4762)(3012,12)
\dottedline{120}(1850,2310)(1850,1812)
\dottedline{120}(1850,2310)(2412,2310) \thicklines
\path(1850,2310)(4812,4046)
\texture{aa777777 77aaaaaa aad5d5d5 d5aaaaaa aa777777 77aaaaaa
aadd5ddd 5daaaaaa
         aa777777 77aaaaaa aad5d5d5 d5aaaaaa aa777777 77aaaaaa aa5ddd5d ddaaaaaa
         aa777777 77aaaaaa aad5d5d5 d5aaaaaa aa777777 77aaaaaa aadd5ddd 5daaaaaa
         aa777777 77aaaaaa aad5d5d5 d5aaaaaa aa777777 77aaaaaa aa5ddd5d ddaaaaaa }

\Thicklines
\dottedline{310}(12,1212)(1850,2310)
\put(1850,2310){\ellipse{60}{60}}

\put(4600,1630){\makebox(0,0)[lb]{\smash{{{\SetFigFont{7}{12.0}{\rmdefault}{\mddefault}{\updefault}$q$}}}}}
\put(2200,4650){\makebox(0,0)[lb]{\smash{{{\SetFigFont{7}{12.0}{\rmdefault}{\mddefault}{\updefault}$\sigma$}}}}}
\put(2200,3130){\makebox(0,0)[lb]{\smash{{{\SetFigFont{7}{12.0}{\rmdefault}{\mddefault}{\updefault}$p$}}}}}

\put(2540,2270){\makebox(0,0)[lb]{\smash{{{\SetFigFont{6}{12.0}{\rmdefault}{\mddefault}{\updefault}$\!\!N$}}}}}

\put(2760,1590){\makebox(0,0)[lb]{\smash{{{\SetFigFont{7}{12.0}{\rmdefault}{\mddefault}{\updefault}$\!p-\!1$}}}}}
\put(1650,1590){\makebox(0,0)[lb]{\smash{{{\SetFigFont{7}{12.0}{\rmdefault}{\mddefault}{\updefault}$-\!1$}}}}}

\put(200,300){\makebox(0,0)[lb]{\smash{{{\SetFigFont{9}{12.0}{\rmdefault}{\mddefault}{\itdefault}Nonexistence}}}}}
\put(0,4300){\makebox(0,0)[lb]{\smash{{{\SetFigFont{9}{12.0}{\rmdefault}{\mddefault}{\itdefault}Existence}}}}}

\put(0,3600){\makebox(0,0)[lb]{\smash{{{\SetFigFont{9}{12.0}{\rmdefault}{\mddefault}{\itdefault}$\boxed{\mu=C_H}$}}}}}
\end{picture}
\vspace{1.5cm}
}

\caption{The nonexistence set $\mathcal{N}$ of equation
\eqref{e:MAIN} for typical values of $\gamma_-$ and $\gamma_+$.}\label{fig}
\end{center}
\end{figure}

\begin{remark}
$(i)$
Observe that in view of the scaling invariance of \eqref{e:MAIN} if
$u(x)$ is a solution to \eqref{e:MAIN} in $B^c_\rho$ then
$\tau^{\frac{p-\sigma}{q-(p-1)}}u(\tau y)$ is a solution to \eqref{e:MAIN}
in $B^c_{\rho/\tau}$, for any $\tau>0$. So in what
follows, for $q\neq p-1$, we confine ourselves to the study of
solutions to \eqref{e:MAIN} on $B^c_1$. For the same
reason, for $q\neq p-1$ we may assume that $C=1$, when convenient.
\smallskip\\
$(ii)$ Using sub- and super-solutions techniques one can show that
if \eqref{e:MAIN} has a positive super-solution in $B_\rho^c$ then
it has a positive solution in $B_\rho^c$ (see Lemma
\ref{l:solSubSuper}). Thus for any
$(q,\sigma)\in\R^2\setminus{\mathcal N}$ equation \eqref{e:MAIN}
admits positive solutions.
\smallskip\\
$(iii)$
Figure \ref{fig} shows the qualitative pictures of the set
${\mathcal N}$ for typical values of $\gamma^-$, $\gamma^+$
and different relations between $p$ and the dimension $N\ge 2$.
\end{remark}


The paper is organized as follows.
Section 2 contains various preliminary results,
including appropriate versions of the Comparison Principle and Weak Maximum Principle in unbounded domains.
In Section 3 we give a new proof of an improved Hardy Inequality with sharp
constants, which is based on Picone's identity and simplifies some arguments
used in the recent papers \cite{AdimChandRamas,Barbatis_Fillipas_Tertikas,Barbatis,GazzolaGrunauMitidieri}.
Section 3 also includes sharp Phragm\'en--Lindel\"of bounds.
The proof of the main result of the paper, Theorem \ref{t:MAIN}, is contained in Section~4.


The Appendix includes various auxiliary results which are systematically used
in the main part of the paper and often are of independent interest.
Parts A and B of the Appendix contain explicit constructions and estimates
of radial sub- and super-solutions to homogeneous $p$-Laplace equations with Hardy--type potentials.
Finally, in Part C of the Appendix we construct {\sl large sub-solutions}
to a homogeneous equation in the exterior of the ball with zero data on the sphere
using the generalized Pr\"ufer transformation technique.

\section{Background, framework and auxiliary facts}

Here and thereafter $N\ge 2$, $1<p<\infty$, $q\in\R$ and $C>0$, unless specified otherwise.
For $0<\rho<R\le +\infty$, we denote the exterior of the closed ball, the open annulus and the sphere
of the radii $\rho$  by
$$B^c_\rho:=\left\{x\in\R^N:\,|x|>\rho\right\},\quad
A_{\rho,R}:=\left\{x\in\R^N:\,\rho<|x|< R\right\},\quad
S_\rho:=\{x\in\R^N:\,|x|=\rho\}.$$
For a function $u=u(x)$ we denote $u^+=\max\{u,0\}$ and $u^-=-\min\{u,0\}$
the positive and negative parts of $u$, respectively.
By $c,c_1,c_2,\dots$ we denote various positive constants whose
exact values are irrelevant.

\paragraph{Homogeneous form associated to $p$-Laplacian.}
Let $\mathcal{E}_V$ be a homogeneous form defined by
\begin{equation}\label{e:FORM}
\mathcal{E}_V(u):=\int_G|\nabla u|^p\,dx-\int_GV|u|^p\,dx
\qquad(u\in W^{1,p}_c(G)\cap C(G)),
\end{equation}
where $G\subseteq\R^N$ is a domain (i.e.\ an open connected set),
and $0\leq V\in L^\infty_{loc}(G)$ a potential. Consider the
equations associated with $\E_V$
\begin{eqnarray}
\label{e:V0}
-\Delta_p u-V|u|^{p-2}u&=&0\quad \mbox{in}\quad G,\\
\label{e:V}
-\Delta_p u-V|u|^{p-2}u&=&f\quad \mbox{in}\quad G,
\end{eqnarray}
where $0\le f\in L^1_{loc}(G)$.
We say that $u\in W^{1,p}_{loc}(G)\cap C(G)$ is a  \textit{super-solution} to equation \eqref{e:V} in a domain
$G\subseteq\R^N$  if for all
$0\leq\varphi\in W^{1,p}_c(G)\cap C(G)$ the following inequality holds
$$\int_G\nabla u|\nabla u|^{p-2}\nabla\varphi
\,dx-\int_G\frac{\mu}{|x|^p}|u|^{p-2}u
\varphi\,dx\geq\int_G f\varphi\, dx.$$ 
The notions of
sub-solution and solution are defined similarly  by replacing
$"\geq"$ with $"\leq"$ and $"="$ respectively.

Let $u\geq 0$ be a solution to \eqref{e:V0} in $G$ and let
$G'\Subset G$. Then the following strong Harnack inequality (cf.
\cite[Theorems 5, 6, 9]{Serrin64}) holds
\begin{equation}\label{pe:sHarnack}
\sup_{G'} u\leq C_S \inf_{G'} u,
\end{equation}
where the constant $C_S>0$ depends on $p,$ $N$, $G'$, $G$  only.
The Harnack inequality and comparison principle in bounded domains \cite{Melian-Sabina de Lis,Takac}
imply that any nontrivial nonnegative super-solution to \eqref{e:MAIN} in $G$
is strictly positive in $G$.

We say that the form $\mathcal{E}_V$ is \textit{positive definite} if
$$\mathcal{E}_V(u)>0,\qquad\forall\:u\in W^{1,p}_{c}(G)\cap C(G), \quad u\neq 0.$$
In this section we study the relation between the positivity of the form $\E_V$ and
the existence of positive super-solutions to the equation \eqref{e:V}.
In the linear case $p=2$ such a  relation is well-documented, see e.g. \cite{Agmon}.
We start with formulating the well--known Picone's Identity for
$p$-Laplacian (see e.g. \cite{AllegrettoHuang,Diaz-Saa,Shafrir}).

\begin{proposition}[\textsf{Picone's Identity}]\label{p:Picon}
Let $w,\phi\in W^{1,p}_{loc}(G)\cap C(G)$ be such that $w\geq 0$ and $\phi>0$. Set
\begin{eqnarray*}
\mathcal{L}(w,\phi)&:=&|\nabla w|^p+(p-1)\left(\frac{w}{\phi}\right)^p
|\nabla\phi|^p-p\left(\frac{w}{\phi}\right)^{p-1}\nabla w|\nabla\phi|^{p-2}\nabla\phi,\\
\mathcal{R}(w,\phi)&:=&|\nabla w|^p-\nabla\left(\frac{w^p}{\phi^{p-1}}\right)
|\nabla\phi|^{p-2}\nabla\phi.
\end{eqnarray*}
Then $\mathcal{L}(w,\phi)=\mathcal{R}(w,\phi)\geq 0$ a.e.\ in $G$.
Moreover, $\mathcal{L}(w,\phi)=\mathcal{R}(w,\phi)=0$ a.e.\ in $G$
if and only if $w=c\phi$ in $G$ for a constant $c>0$.
\end{proposition}

An immediate consequence of Picone's identity is that the existence of a positive super-solution
to \eqref{e:V} implies positivity of the form $\E_V$, as the following proposition shows.

\begin{proposition}\label{p:form+}
Let $\phi>0$ be a super-solution (sub-solution) to equation \eqref{e:V}.
Then the form $\E_V$ satisfies the following inequality
\begin{equation}\label{e:FormRepr}
\mathcal{E}_V(u)\geq(\leq)\,\int_G
\mathcal{R}(u,\phi)\,dx+\int_G\frac{f}{\phi^{p-1}}|u|^p\,dx,
\qquad\forall\:u\in W^{1,p}_c(G)\cap C(G).
\end{equation}
\end{proposition}
\begin{remark}\label{r:E>0}
$({i})$
If $\phi>0$ is a super-solution to \eqref{e:V} then
$\mathcal{R}(u,\phi)\ge 0$ a.e.\ in $G$ and, in particular,
$$\mathcal{E}_V(u)\geq\int_G \frac{f}{\phi^{p-1}}|u|^p\,dx\geq0,
\qquad\forall\:u\in W^{1,p}_c(G)\cap C(G).$$ If, in addition,
$f>0$ then
$$\mathcal{E}_V(u)>0,
\qquad\forall\:u\in W^{1,p}_c(G)\cap C(G),\quad u\neq 0.$$
$({ii})$
If $\phi>0$ is a solution to \eqref{e:V} then inequality
\eqref{e:FormRepr} becomes an identity.
\end{remark}
\begin{proof}
Let $\phi>0$ be a super-solution (sub-solution) to \eqref{e:V}.
Testing \eqref{e:V} by $\xi=\frac{|u|^p}{\phi^{p-1}}\in W^{1,p}_c(G)\cap C(G)$ we obtain
$$\int_G Vu^p\, dx\,\leq(\geq)\,
p\,\int_G\frac{u\nabla\phi}{\phi}\,\left|\frac{u\nabla\phi}{\phi}\right|^{p-2}\nabla u\,dx-
(p-1)\int_G|\nabla\phi|^p \left|\frac{u}{\phi}\right|^{p}dx-
\int_G \frac{f}{\phi^{p-1}}|u|^p\,dx,$$
which implies \eqref{e:FormRepr}.
\end{proof}

The following straightforward corollary of Proposition \ref{p:form+}
is our main tool in proving nonexistence of positive solutions
to nonlinear equation \eqref{e:MAIN}.

\begin{corollary}[\textsf{Nonexistence principle}]\label{c:form-}
Assume that there exists $u\in W^{1,p}_c(G)\cap C(G)$ such that
$\mathcal{E}_V(u)<0$.
Then equation \eqref{e:V0} has no positive super-solution.
\end{corollary}

Another interesting application of Proposition \ref{p:form+}
is a version of Barta's inequality (cf. \cite{AllegrettoHuang}).
\begin{corollary}[\textsf{Barta's inequality}]\label{c:Barta's}
Assume that equation \eqref{e:V0} admits a positive
super-solu\-tion. Then for every  $0<\varphi\in W^{1,p}_{loc}(G)\cap C(G)$
such that $-\Delta_p \varphi-V\varphi^{p-1}\in L^1_{loc}(G)$ the
following inequality holds
\begin{equation}\label{e:34}
\inf_{x\in G}\frac{-\Delta_p \varphi-V\varphi^{p-1}}{\varphi^{p-1}}\le
\inf_{0\lneq u\in W^{1,p}_{c}\cap C(G)}\frac{\int_G(|\nabla u|^p-V u^p)dx}{\int_G u^p\,dx}.
\end{equation}
\end{corollary}
\begin{proof}
Set $F(x):=-\Delta_p \varphi-V\varphi^{p-1}\in  L^1_{loc}(G)$. We
may assume that $F\ge 0$ (otherwise inequality \eqref{e:34} is
trivial). Proposition \ref{p:form+} implies that
$$\E_V(u)\geq
\int_G\frac{F}{\varphi^{p-1}} u^p dx \geq \inf_{x\in G}\frac{F}{\varphi^{p-1}}\int_G u^p dx,
\qquad\forall\:u\in W^{1,p}_c(G)\cap C(G).$$
So the assertion follows.
\end{proof}

We need the following version of the Caccioppoli inequality, which is
a consequence of Proposition \ref{p:form+}.

\begin{corollary}[\textsf{Caccioppoli-type Inequality}]\label{c:Caccioppoli}
Let $u>0$ be a sub-solution to \eqref{e:V0}. Then
\begin{equation}\label{e:Caccioppoli}
\int_{G}|\theta\,\nabla u|^pdx\leq p\int_G Vu^p|\theta|^pdx+p^p\int_G u^p|\nabla\theta|^pdx,
\qquad\forall\:\theta\in W^{1,\infty}_c(G).
\end{equation}
\end{corollary}

\begin{proof} From \eqref{e:FormRepr} we have
\begin{eqnarray*}\mathcal{E}_V(u\theta)&\leq& \int_G|\nabla (u\theta)|^pdx-
p\int_G\theta \nabla u \,|\theta \nabla u|^{p-2}\nabla (u\theta)\,
dx+ (p-1)\int_G|\theta\nabla u|^p dx \\
&\leq&
\int_G|\nabla (u\theta)|^pdx+p\int_G |\theta\,\nabla u|^{p-1} u|\nabla\theta|\,dx-\int_G|\theta\,\nabla u|^p dx.
\end{eqnarray*}
Using the Young's Inequality and \eqref{e:FORM} we obtain
$$\int_G|\theta\nabla u|^pdx\leq\int_G V u^p|\theta|^pdx+ p^{p-1}\int_G|u\,\nabla\theta|^pdx
+\frac{p-1}{p}\int_G |\theta\,\nabla u |^pdx,$$
so the assertion follows.
\end{proof}

\paragraph{Comparison and Maximum Principles.}
We say that $0\leq w\in W^{1,p}_{loc}(G)$ satisfies condition $(\mathcal{S})$ if the following holds:
\begin{enumerate}\label{e:Shafrir}
\item[$(\mathcal{S})$] there exists $(\theta_n)_{n\in \N}\subset
W_c^{1,\infty}(\R^N)$ such that $0\leq\theta_n\to 1$ a.e. in $\R^N$ and
$$\int_{G}\mathcal{R}(\theta_nw,w)\,dx \to 0\quad\mbox{as }\:n\to +\infty.$$
\end{enumerate}
Notice that if $G$ is bounded and $w\in W^{1,p}(G)$ then  condition $(\mathcal{S})$ is
trivially satisfied  with $\theta=1$ in $G$.

Using condition $(\mathcal S)$,
we establish a version of comparison principle in a form suitable for our framework.
The proof follows with certain modifications the ideas in \cite{ShafrirMarcus,ShafrirPoljakovsky,Shafrir}.

\begin{theorem}[\textsf{Comparison Principle}]\label{l:WCPimpr}
Let $q<p-1$ and $0\leq f\in L^1_{loc}(G)$.
Let $0<u\in W^{1,p}_{loc}(G)\cap C(\bar G)$ be a super-solution
and $v\in W^{1,p}_{loc}(G)\cap C(\bar G)$ a sub-solution to
equation
\begin{eqnarray}\label{e:V'}
-\Delta_p u-V|u|^{p-2}u&=&f|u|^{q-1}u\quad \mbox{in}\quad G.
\end{eqnarray} If $G$ is an unbounded domain, assume in
addition that $\partial G\neq \emptyset$ and $v^+$ satisfies condition $(\mathcal{S})$.
Then $u\ge v$ on $\partial G$ implies $u\ge v$ in $G$.
\end{theorem}
\begin{proof}
Let $S:=\{x\in G\,:\,v>u\}$.
Assume  for a contradiction that $S\neq\emptyset$.
Then
$$K:=\sup_{x\in S}\left(\log \frac{v}{u}\right)\in(0,+\infty].$$
Fix a positive constant $b$ such that $5b<K$.
Let $\eta\in C^1(\R)$ be a nondecreasing function  such that
$$\text{$\eta(t)=0$ for $t\leq 2b$,\quad $\eta(t)=1$ for $t\geq 5b$\quad and
\quad$\eta^\prime(t)>0$ for $3b\leq t\leq 4b$.}$$ Let
$\xi=\eta(\log\frac{v}{u})$. Then $0\le \xi\in
W^{1,p}_{loc}(G)\cap C(\bar G)$ and $\spt(\xi)\subset S\subseteq
G$. Let $\theta\in W^{1,\infty}_c(\R^N)$ 
then $\spt(\theta\,\xi)$ is  compact  in $G$. Later on we specify $\theta$ for
the case of  a bounded and unbounded $G$.
Set
$$\phi_1:=\left(\frac{\theta^pv^p}{u^{p-1}}\right)\xi,\qquad
\phi_2:=\theta^pv\xi.$$
Clearly $\phi_1,\,\phi_2\in W^{1,p}_{c}(G)$.
Since $u$ is a super-solution to \eqref{e:V'}, testing
\eqref{e:V'} by $\phi_1$ and using Picone's Identity  we infer that
\begin{eqnarray}\label{e:51}
0&\leq &
\int_S \xi|\nabla u|^{p-2}\cdot\nabla
u\cdot\nabla\left(\frac{\theta^pv^p}{u^{p-1}}\right)\,dx+
\int_S \theta^pv^p\,\nabla \log u\cdot|\nabla\log u |^{p-2}\cdot\nabla\xi\,dx\nonumber\\
&-&
\int_S V \theta^pv^p\xi dx-\int_S\frac{fu^q}{u^{p-1}}\theta^pv^p\xi\,dx\nonumber\\
&=& \int_S|\nabla(\theta v)|^p\xi\,dx-\int_S \mathcal{R}(\theta v,u)\xi\,dx
+\int_S \theta^pv^p\,\nabla \log u\cdot|\nabla\log u |^{p-2}\cdot\nabla\xi\,dx\nonumber\\
&-&\int_S V \theta^pv^p\xi dx -\int_S
fu^{q-(p-1)}{\theta^pv^p}\xi\,dx.\nonumber
\end{eqnarray}
Thus from Proposition \ref{p:Picon} we obtain
\begin{equation*}
\int_S|\nabla(\theta v)|^p\xi\,dx+ \int_S \theta^pv^p|\nabla\log u |^{p-2}
\nabla \log u\cdot\nabla\xi\,dx-
\int_S V \theta^pv^p\xi\,dx\ge\int_S {f}{u^{q-(p-1)}}\theta^pv^p\xi\,dx.
\end{equation*}
Since $v$ is a sub-solution to \eqref{e:V'}, testing \eqref{e:V'} by $\phi_2$ we derive
\begin{equation*}
\int_S \xi|\nabla v|^{p-2}\nabla v \cdot\nabla (\theta^p v)dx
+\int_S \theta^pv^p|\nabla\log v|^{p-2}\nabla\log v\cdot\nabla\xi dx
-\int_S V\theta^pv^p\xi dx
\le\int_Sf\theta^pv^{q+1}\xi dx.
\end{equation*}
Subtracting the former inequality 
from the latter  one and using Picone's Identity again 
we obtain
\begin{eqnarray}\label{e:16}
&&\mathcal I_\theta:= \int_S\theta^pv^p\,\left(|\nabla
\log v |^{p-2}\cdot\nabla\log v
-|\nabla \log u |^{p-2}\cdot\nabla\log u\right)\nabla\xi\,dx\\
&&\leq\int_S \left\{|\nabla(\theta v)|^p -|\nabla
v|^{p-2}\cdot\nabla v\cdot\nabla (\theta^pv)\right\}\xi\,dx-
\int_S {f}\theta^pv^p(u^{q-(p-1)}-v^{q-(p-1)})\xi\,dx\nonumber\\
&&\le\int_{S\cap \spt(\theta\xi)} \mathcal{R}(\theta v,v)\,dx.\nonumber
\end{eqnarray}
We claim that
\begin{equation}\label{e-claim}
\mathcal I_\ast:=\int_S v^p\left(|\nabla\log
v|^{p-2}\cdot\nabla\log v-|\nabla\log u|^{p-2}\cdot\nabla\log
u\right)\nabla\xi\,dx\le 0
\end{equation}
implies $S=\emptyset$. Define the open subset $S'\subset S$ by
 $$S^\prime:=\{x\in G:\,\left(\log \frac{v}{u}\right)\in (3b,4b)\}\subset S$$
and observe that $\eta'(\log\frac{v}{u})>0$ on $S'$. There exists at least one connected component  $S_{i}$ of the set $S'$ such that
$\log\frac{u}{v}$ attains  all values between $3b$ and $4b$ on $S_{i}$.

Since
$$\nabla\xi=(\nabla \log v-\nabla\log u) \,\eta'(\log\frac{v}{u}),$$
and
$$\left(|z_1|^{p-2}z_1-|z_2|^{p-2}z_2\right)(z_1-z_2)\geq 0,\qquad\forall\,z_1,z_2\in\R^N,$$
with equality if and only if $z_1=z_2$, from \eqref{e-claim} we have $\mathcal I_\ast=0$.
Therefore  $\log\frac{v}{u}=c_{i}$ 
on $S_{i}$, which is a contradiction.
\smallskip

Below we show that \eqref{e-claim} holds. Indeed, if the domain $G$ is bounded then
$\spt(\xi)\Subset S$ and one simply chooses $\theta\equiv 1$
on $\bar G$ and $\theta\equiv 0$ on $R^N\setminus\bar G$. Then
\eqref{e:16} implies that $\mathcal I_\ast\le 0$.

Now let $G$ be  an unbounded domain. Let $\theta_n$ satisfies
condition $(\mathcal{S})$. Then $\spt(\theta_n)\cap\spt(\xi)\neq
\emptyset$  for $n$ large enough and from \eqref{e:16} we have
$$I_{\theta_n}\leq \int_G\mathcal{R}(\theta_n v,v)\,dx\to 0\qquad\text{as}\quad n\to\infty.$$
So the assertion follows.
\end{proof}

The proof of the following lemma follows closely the arguments in \cite[Lemma 2.9]{Agmon-2}.

\begin{lemma}\label{p:u^-_sub}
Let $v$ be a sub-solution to \eqref{e:V0}.
Then $v^+$ is a sub-solution to \eqref{e:V0}.
\end{lemma}
\begin{proof}
For any $\epsilon>0$ define
$v_\epsilon=(v^2+\epsilon^2)^{{1}/{2}}$.
Then $0< v_\epsilon\in W^{1,p}_{loc}(G)$ and
by the Lebesgue dominated convergence theorem,  $v_\epsilon$ converges to $|v|$ in $W^{1,p}_{loc}(G)$.
Let $0\leq\phi\in W^{1,p}_{loc}(G)\cap C(G)$. A direct computation shows that
$$\nabla v_\epsilon\cdot\nabla \phi=\nabla v\cdot\nabla \left(\frac{v}{v_\epsilon}\,\phi\right)-
\frac{v_\epsilon^2-v^2}{v_\epsilon^3}\,\phi\,|\nabla v|^2,$$ which
implies that
$
\nabla v_\epsilon\cdot\nabla \phi\leq\nabla v\cdot\nabla
\left(\frac{v}{v_\epsilon}\,\phi\right)$.
Set
$\phi_\epsilon=\frac{1}{2}\left(1+\frac{v}{v_\epsilon}\right)\phi.$
It follows that
\begin{equation}\label{e:19}
\frac{1}{2}\,\nabla(v+v_\epsilon)\cdot\nabla\phi=\frac{1}{2}\left(\nabla
v\cdot\nabla\phi+ \nabla v_\epsilon\cdot\nabla\phi\right)\leq
\nabla v\cdot\nabla\left( \frac{1}{2}\left(1+ \frac{
v}{v_\epsilon}\right)\phi\right)=\nabla
v\cdot\nabla\phi_\epsilon.\end{equation}
Testing \eqref{e:V0} against $\phi_\epsilon$ and using \eqref{e:19} we derive
\begin{eqnarray}\label{e:20}0&\geq&\int_G|\nabla
v|^{p-2}\nabla v\cdot\nabla\phi_\epsilon \,dx-\int_G V
|v|^{p-2}v\phi_\epsilon\,dx\\
&\geq&\int_G|\nabla v|^{p-2}\,\nabla
\frac{1}{2}\,(v+v_\epsilon)\cdot\nabla\phi \,dx-\int_G V
|v|^{p-2}v\phi_\epsilon\,dx.\nonumber
\end{eqnarray}
Notice that $v+v_\epsilon\rightarrow v^+ $ and
$\phi_\epsilon\rightarrow \chi_{\{v^+>0\}}\phi$  a.e.\ in $G$ as
$\epsilon\rightarrow 0$. Letting $\epsilon\rightarrow 0$
in \eqref{e:20} we infer that
$$\int_G|\nabla
(v^+)|^{p-2}\,\nabla (v^+)\cdot\nabla\phi \,dx-\int_G V(
v^+)^{p-1}\phi\,dx\leq 0,$$
which completes the proof.
\end{proof}

We establish the Weak Maximum Principle for super-solutions to \eqref{e:V0}
as a corollary of the Comparison Principle and Lemma \ref{p:u^-_sub}.

\begin{proposition}[\textsf{Weak Maximum Principle}]
Let $\partial G\neq \emptyset$. Assume that \eqref{e:V0}  admits a positive
super-solution $0<\phi\in W^{1,p}_{loc}(G)\cap C(\bar G)$.
Let $u\in W^{1,p}_{loc}(G)\cap C(\bar G)$ be a super-solution to
equation \eqref{e:V0} such that $u\ge 0$ on $\partial G$. For an
unbounded $G$ assume in addition that $u^-$ satisfies
condition $(\mathcal{S})$. Then $u\geq
0$ in $G$.
\end{proposition}
\begin{proof}
By Proposition \ref{p:u^-_sub} observe that $u^-\in
W^{1,p}_{loc}(G)\cap C(\bar G)$ is a sub-solution to \eqref{e:V0} and
$u^-=0$ on $\partial G$. Thus $u^-\le\eps\phi$ on $\partial G$,
for any $\eps>0$. By Theorem \ref{l:WCPimpr}, we conclude that
$u^-\le\eps\phi$ in $G$ for an arbitrary small $\eps>0$. Hence
$u^-=0$ in $G$.
\end{proof}

\begin{remark}
After examples constructed in \cite{DelPino,Fleck-1}
(see also discussions in \cite{Melian-Sabina de Lis,Takac})
it is known that the form $\E_V$ is nonconvex as soon as $p\neq 2$ and the potential $V$
has a nontrivial 'negative' part $V^+$, even if $\E_V$ is nonnegative and admits
representation \eqref{e:FormRepr} with respect to a positive super-solution of \eqref{e:V}.
One of consequences of this fact is that the assumption $f\ge 0$ in Theorem \ref{l:WCPimpr}
can not be removed, otherwise the comparison principle fails.
\end{remark}

\paragraph{Positive solution between sub- and super-solutions.}
We show that the existence of a positive super-solution to nonlinear equation \eqref{e:MAIN}
implies the existence of a positive {\sl solution} to \eqref{e:MAIN}.
The following result on bounded domains is standard.

\begin{lemma}\label{l:WCP_c}
Let $\mu\leq C_H$ and $G\subset B_1^c$ be a bounded smooth domain.
Let $v,u\in W^{1,p}(G)\cap C(\bar G)$ be a sub- and  super-solution to \eqref{e:MAIN} in $G$, respectively.
Assume that $0<v\le u$ in $\bar G$.
Then there exists a solution $w\in W^{1,p}(G)\cap C(\bar G)$ to \eqref{e:MAIN} in $G$,
so that $v\leq w\leq u$ in $G$ and $w=v$ on $\partial G$.
\end{lemma}
\begin{proof}
The proof is a standard consequence of the comparison principle and monotone iterations scheme
(cf. \cite{Drabek,Takac} for similar results). We omit the details.
\end{proof}

By means of the standard digitalization techniques Lemma \ref{l:WCP_c} extends to the following.

\begin{proposition}\label{l:solSubSuper}
Let $\mu\leq C_H$.
Assume that \eqref{e:MAIN} has a positive super-solution in $B^c_1$.
Then \eqref{e:MAIN} has a positive solution in $B^c_1$.
\end{proposition}
\begin{proof} Let $u>0$ be a super-solution to \eqref{e:MAIN}.
Set $v=cr^{\gamma_-}$ and observe that
$$-\Delta_p v-\frac{\mu}{|x|^p}{v}^{p-1}=0\quad \mbox{in}\quad B^c_{1},$$
so $v>0$ is a sub-solution to nonlinear equation \eqref{e:MAIN} in $B^c_1$.
By Proposition \ref{p:R(v,v)}, $v$ satisfies  condition $(\mathcal{S})$.
Choose $c$ in such a way that $u\geq cv$ for $|x|=2$.
Thus Theorem \ref{l:WCPimpr} implies that $u\geq v$ in $B^c_{2}$.
By Lemma \ref{l:WCP_c}, for each $n\ge 3$ there exists a solution
$w_n\in W^{1,p}(A_{2,n})$ to \eqref{e:MAIN} in $A_{2,n}$ such that
\begin{eqnarray}\label{e:71}
v\leq w_n\leq u\quad\mbox{in}\quad A_{2,n},
\qquad w_n=v\quad \mbox{on}\quad\partial A_{2,n}.
\end{eqnarray}
By Corollary \ref{c:Caccioppoli},
we conclude that there exists a constant $M_n>0$ such that
\begin{equation}\label{e:00}
||\nabla w_{n+1}||_{L^p(A_{3,n})}\leq M_{n},\qquad\forall\, n\geq 4.
\end{equation}
Using \eqref{e:00} and \eqref{e:71},
one can proceed following the standard digitalization techniques
in order to construct a solution to \eqref{e:MAIN} with the required properties.
\end{proof}

\section{Hardy inequalities and positive super-solutions}\label{s-Hardy}

One of the crucial components in our proof of Theorem \ref{t:MAIN}
is an improved Hardy inequality on exterior domains.
Inequalities of this type were recently obtained by several authors
using various techniques,
see \cite{Peral-2,AdimChandRamas,Barbatis_Fillipas_Tertikas,Barbatis,GazzolaGrunauMitidieri}.
Here we give a simple proof of an improved Hardy inequality on exterior domains
for all $p>1$ and $N \geq 2$, which is based on the explicit construction of
appropriate super- and sub-solution and transformation \eqref{e:FormRepr}.

Throughout the paper we use the notation $\gamma_*:=\frac{p-N}{p}$ and
\begin{equation}\label{E:HardyNotations}
C_H:=\left|\frac{p-N}{p}\right|^p,
\quad
C_*:=
\left\{\begin{array}{ll}
\frac{p-1}{2p}\left|\frac{N-p}{p}\right|^{p-2},&\text{
$N\neq p$},\smallskip\\
\left(\frac{N-1}{N}\right)^N,&\text{
$N=p$},\\
\end{array}
\right.\quad
m_*:=\left\{
\begin{array}{ll}
2,&\text{
$N\neq p$},\\
N,&\text{
$N=p$}.\\
\end{array}
\right.
\end{equation}
Recall that, according to Proposition \ref{p:form+}
the existence of a positive super-solution to the equation
\begin{equation}\label{e:Hardy}
-\Delta_pu-\frac{\mu}{|x|^p}u^{p-1}-\frac{\epsilon}{|x|^p\log^{m_*}|x|}u^{p-1}=0\quad
\mbox{in}\quad B^c_\rho,
\end{equation}
with some $\rho>1$, implies that the form
\begin{equation*}
\mathcal{E}_{\mu,\eps}(u)=\int_G|\nabla u|^pdx-
\mu\int_G\frac{u^p}{|x|^p}dx-\epsilon\int_G\frac{|u|^p}{|x|^p\log^{m_*}|x|}dx
\qquad (u\in W^{1,p}_c(B^c_\rho)\cap C(B^c_\rho)),
\end{equation*}
is nonnegative.
Thus, in order to prove an improved  Hardy inequality it is sufficient
to find a super-solution for the corresponding equation.
The idea to use Picone's identity for proving Hardy type inequalities related to
$p$-Laplace operator goes back to \cite{AllegrettoHuang}, see also \cite{Peral-1,Peral-2}. 
However, as discovered in \cite{Harrell}, 
such a technique can be in fact attributed as far as to an 1907's paper by Boggio \cite{Boggio}.

\begin{theorem}[\textsf{Improved Hardy Inequality}]\label{l:HardyImpr}
For every $p>1$ there exists $\rho\geq 1$ such that
\begin{equation}\label{e:HardyImpr}\int_{B^c_\rho}|\nabla v|^pdx\geq
C_H\int_{B^c_\rho}\frac{|v|^p}{|x|^p}dx+
C_*\int_{B^c_\rho}\frac{|v|^p}{|x|^p\log^{m_*}|x|}dx,
\qquad\forall\:v\in W^{1,p}_c(B^c_\rho)\cap C(B^c_\rho).
\end{equation}
The constants $C_H$ and $C_*$ are sharp in the sense that the
inequality
$$\E_{\mu,\epsilon}(v)\ge 0,
\qquad\forall\:v\in W^{1,p}_c(B^c_\rho)\cap C(B^c_\rho),$$
fails in any of the following two cases:
\begin{enumerate}
\item[{$(i)$}] $\mu=C_H$, $\epsilon>C_*$,
\item[{$(ii)$}] $\mu>C_H$, $\epsilon\in\R$.
\end{enumerate}
\end{theorem}

\begin{proof}
Lemma \ref{p:SubSuper'}  for $p\neq N$ and a direct
computation for $p=N$ show that the function
$$\phi(r)=r^{\gamma_*}(\log r)^\beta (\log\log r)^\tau,
\qquad
\left\{
\begin{array}{ccl}
\beta=\frac{1}{p},&\tau\in(0,\frac{2}{p}) & \text{for $p\neq N$},\smallskip\\
\beta=\frac{N-1}{N},&\tau=0 & \text{for $p=N$}.
\end{array}
\right.
$$
is a super-solution to equation \eqref{e:Hardy} with $\mu=C_H$ and
$\epsilon=C_\ast$ in $B^c_\rho$ with some $\rho> 1$. Thus
\eqref{e:HardyImpr} follows immediately from Proposition
\ref{p:form+}.

\paragraph{{\sl Sharpness of the constants.}}(i)
Define
$$\phi(r)=r^{\gamma_\ast}(\log r)^\beta (\log\log r)^\tau,
\qquad
\left\{
\begin{array}{lcl}
\beta=\frac{1}{p},&\tau\in(-\frac{1}{p},0) & \text{for $p\neq N$},\smallskip\\
\beta=\frac{N-1}{N},&\tau=0 & \text{for $p=N$}.
\end{array}
\right.
$$
By Lemma \ref{p:SubSuper'} (ii) one can choose $\rho>1 1$
such that $\phi$ is a sub-solution with $\mu=C_H$ and
$\epsilon=C_\ast$ in $B^c_\rho$. Let $R>\rho$. Following
\cite{AdimChandRamas}, we define the  cut-off function
\begin{equation}\label{e:theta_log}
\theta_{R}(t):=\left\{
\begin{array}{cl}
2t/\rho-3, & \frac{3}{2}\rho\leq t\leq 2\rho,\smallskip\\
1 & 2\rho\leq t\leq R,\\
\frac{\log\frac{R^2}{t}}{\log R},& R\leq t\leq R^2.
\end{array}\right.
\end{equation}
Below we show that for any $\eps>0$,
$$\E_{C_H,C_\ast+\eps}(\phi\,\theta_{R}^\alpha)\to-\infty\quad\text{as}\quad R\to\infty,$$
where $\alpha=1$ if $p\geq 2$, and $\alpha>\frac{2}{p}$ if $p<2$.
By Proposition \ref{p:form+} and using
\eqref{B1-ii}, \eqref{B1-ii-2}, \eqref{B1-iii} we obtain
\begin{eqnarray*}
\E_{C_H,C_\ast}(\phi\,\theta_R^\alpha)&\leq&
c_1+c_2\int_{A_{R,R^2}}
\mathcal{R}(\theta_R^\alpha\,\phi,\phi)\,dx\leq c_3.
\end{eqnarray*}
Further, it is easy to see that
\begin{eqnarray*}
\int_{A_{\frac{3}{2}\rho,R^2}}\frac{|\phi\,\theta_R^\alpha|^p}{|x|^p\log^{m_\ast}|x|}dx&\geq&
\int_{2\rho}^R \frac{(\log\log r)^{\tau p}}{r\log r }dr=
c_4(\log\log R)^{\tau p+1}-c_5.
\end{eqnarray*}
Thus for any $\eps>0$ we arrive at
\begin{eqnarray*}
\E_{C_H,C_\ast+\eps}(\phi\,\theta_R^\alpha)=
\E_{C_H,C_\ast}(\phi\,\theta_R^\alpha)-\eps\int_{B^c_{\frac{3}{2}\rho}}\frac{|\phi\,\theta_R^\alpha|^p}{|x|^p\log^2|x|}dx
\to -\infty \quad\mbox{as}\quad R\to\infty.
\end{eqnarray*}

(ii)
Choosing $\phi(r)=r^{\gamma_\ast}$ as a sub-solution
to \eqref{e:Hardy} with $\mu=C_H$ and $\epsilon=C_\ast$ in $B^c_2$,
one can verify that \eqref{e:HardyImpr} with $\mu>C_H$ and any $\epsilon\in \R$
fails on the family of functions $\phi\theta_{R}$ defined as above.
\end{proof}

As a consequence of the last theorem we obtain
the following nonexistence result,
which is crucial in our proofs of nonexistence of positive super-solutions
to nonlinear equation \eqref{e:MAIN}.
\begin{corollary}\label{c:NonExist}
Equation \eqref{e:Hardy} admits positive super-solutions in $B^c_\rho$ with some $\rho>1$
if and only if $\mu<C_H$ and $\epsilon\in\R$, or $\mu=C_H$ and $\epsilon\le C_*$.
\end{corollary}

\begin{remark}
Equation \eqref{e:Hardy} with $\epsilon\neq 0$ is not homogeneous with respect to scaling,
 i.e.\ the existence of a positive
(super)\,solution in $B^c_\rho$ with $\rho>1$ does not imply the existence
of positive (super)\,solution in $B^c_1$ and so
the value of the radius $\rho>1$ becomes essential.
\end{remark}

Next we describe the behavior at infinity of positive super-solutions to equation \eqref{e:Hardy}
in the case when $\mu\leq C_H$ and $\epsilon\in[0,C_*)$.
For $\epsilon\in[0,C_*)$, denote by $\beta_-<\beta_+$ the real roots of the equation
\begin{equation}\label{e:1000}
\begin{array}{lccl}
\frac{1}{2}|\gamma_\ast|^{p-2}(p-1)(2-\beta p)\beta&=&\epsilon&\mbox{if}\quad p\neq N,\smallskip\\
(N-1)(1-\beta)\beta^{N-1}&=&\epsilon&\mbox{if}\quad p=N.
\end{array}
\end{equation}
Notice that $0\le\beta_-<\frac{1}{p}<\beta_+\le\frac{2}{p}$ if
$p\neq N$ and $0\le\beta_-<\frac{N-1}{N}<\beta_+\le 1$ if $p=N$.

\begin{theorem}[\textsf{Lower bound}]\label{l:WCP}
Let $u>0$ be a super-solution to \eqref{e:Hardy} in $B^c_\rho$.
The following assertions are valid.
\begin{enumerate}
\item[{$(i)$}]
Let $\mu\leq C_H$, $\epsilon=0$.  There exists $c>0$ such that
$$u\geq c |x|^{\gamma_-},\quad x\in B^c_{2\rho}.$$
\item[{$(ii)$}]
Let $\mu=C_H$, $\epsilon=0$.  There exists $c>0$ such that
$$u\geq c |x|^{\gamma_\ast},\quad x\in B^c_{2\rho}.$$
\item[{$(iii)$}]
Let $\mu=C_H$, $\epsilon\in(0,C_*)$. For every $\tau<0$ there exists $c>0$ such that
$$u\geq c |x|^{\gamma_\ast}(\log |x|)^{\beta_-}(\log\log |x|)^\tau,\quad x\in B^c_{2\rho}.
$$
\end{enumerate}
\end{theorem}
\begin{proof}
Follows from Theorem \ref{l:WCPimpr} and
small sub-solutions estimates in Proposition \ref{p:R(v,v)}.
\end{proof}

The next lemma establishes a Phragm\'en--Lindel\"of type
upper bound on super-solutions.

\begin{theorem}[\textsf{Upper bound}]\label{l:upper_est<C_H}
Let $u>0$ be a super-solution to \eqref{e:Hardy} in $B^c_\rho$. The following
assertions are valid.
\begin{enumerate}
\item[{$(i)$}]
Let $\mu<C_H$, $\epsilon=0$.
There exists $c>0$ such that
$$\inf_{S_R}u\leq c R^{\gamma_+},\quad R>2\rho.$$
\item[{$(ii)$}]
Let $\mu=C_H$, $\epsilon=0$.
There exists $c>0$ such that
$$\inf_{S_{R}}u\leq c R^{\gamma_\ast}(\log R)^{\beta^*},
\quad R>2\rho,$$
where
$\beta^*=\frac{2}{p}$ for $p\neq N$ or $\beta^*=1$ for $p=N$.
\item[{$(iii)$}] Let  $\mu=C_H, 0<\epsilon<C_*$.
For every $\beta\in(\beta_+,\beta_*)$ there exists $c>0$ such that
$$\inf_{S_{R}}u\leq cR^{\gamma_\ast} (\log R)^{\beta},
\quad R>2\rho.$$
\end{enumerate}
\end{theorem}
\begin{proof}
Let $v>0$ be a large sub-solution to \eqref{e:Hardy},
that is a positive sub-solution to \eqref{e:Hardy}
that satisfies the boundary condition $v=0$ on $S_\rho$,
as constructed in Appendix \ref{A:B}.
We are going to show that
\begin{equation}\label{e:62}
\inf_{S_R}u\leq c\sup_{S_R} v,
\quad R> 2\rho.
\end{equation}
For a contradiction, assume that for an arbitrary large $c>0$ there exists $R>2\rho$
so that $u\geq cv$ on $S_R$. Thus
$$u-cv\geq 0\quad\mbox{on }\:\partial A_{\rho,R}.$$
Then Theorem \ref{l:WCPimpr}, applied on $A_{\rho,R}$ yields
$$u-cv\geq 0\quad\mbox{on }\: A_{\rho,R}.$$
In particular, this implies that
$$ u(x)\geq c v(x),\quad x\in S_{2\rho}.$$
But this contradicts to the continuity of $u$.

Now the assertions (i)-(iii) follow from \eqref{e:62} via Theorems \ref{t:B0} and \ref{t:B1}.
\end{proof}

\section{Proof  Theorem \ref{t:MAIN}}

First, we prove the nonexistence of positive super-solutions to \eqref{e:MAIN}
in the super-homogeneous case $q\geq p-1$ and sub-homogeneous case $q<p-1$.
After this we show sharpness of our nonexistence results by
constructing explicit super-solutions in all complementary cases.

\subsection{Nonexistence: super-homogeneous case $q\geq p-1$}
We distinguish between the cases $\mu<C_H$ and $\mu=C_H$.
\paragraph{Case $\mu<C_H$.}
First we prove the  nonexistence of super-solutions in the subcritical case,
i.e.\ when $(q,\sigma)$ is below the critical line $\Lambda^*$.

\begin{proposition}
Let  $\sigma<\gamma_-(q-p+1)+p$.
Then \eqref{e:MAIN} has no positive super-solution in $B^c_1$.
\end{proposition}
\begin{proof}
Let $u>0$ be a super-solution to \eqref{e:MAIN} in $B^c_1$.
Then $u$ is a super-solution to the homogeneous equation
\begin{equation}\label{e:01}
-\Delta_p u-\frac{\mu}{|x|^p} u^{p-1}= 0\quad \mbox{in}\quad B^c_1.
\end{equation}
By Theorem \ref{l:WCP}(i) we conclude that $u\geq c_1|x|^{\gamma_-}$ in $B^c_2$.
Thus from equation \eqref{e:MAIN} it follows that $u>0$ is a super-solution to
\begin{equation} \label{e:02}
-\Delta_p u -\frac{\mu+W(x)}{|x|^p}u^{p-1}=0\quad \mbox{in}\quad B^c_2,
\end{equation}
where
$$W(x):=C|x|^{p-\sigma}u^{q-(p-1)}\geq Cc_1^{q-(p-1)}|x|^{{\gamma_-}(q-p+1)+p-\sigma},$$
with $\gamma_-(q-p+1)+p-\sigma>0$.
Then the assertion follows by Corollary \ref{c:NonExist}.
\end{proof}

Next we prove the nonexistence  in the critical case,
i.e.\ when $(q,\sigma)$ belongs to the critical line $\Lambda^*$.

\begin{proposition}
Let  $\sigma=\gamma_-(q-p+1)+p$.
Then \eqref{e:MAIN} has no positive super-solution in $B^c_1$.
\end{proposition}
\begin{proof}
Let $u>0$ be a super-solution to \eqref{e:MAIN} in $B^c_1$.
Arguing as in the proof above,
we conclude that $u$ is a super-solution to \eqref{e:02}, where
$$W(x):=C|x|^{p-\sigma}u^{q-(p-1)}\geq Cc_1^{q-(p-1)}|x|^{{\gamma_-}(q-p+1)+p-\sigma}=c.$$
Thus $u>0$ is a super-solution to the homogeneous equation
\begin{eqnarray*}
-\Delta_p u-\frac{\tilde{\mu}}{|x|^p}u^{p-1}=0\quad \mbox{in}\quad B^c_2,
\end{eqnarray*}
where $\tilde{\mu}=\mu+c$.
Without loss of generality, we may assume that $\tilde{\mu}<C_H$.
Then by Theorem~\ref{l:WCP}(i) we conclude that $u\geq c_2|x|^{\tilde\gamma_-}$ in $B^c_2$,
with $\tilde\gamma_-\in(\gamma_-,\gamma_\ast)$.
Therefore $u$ is a super-solution to~\eqref{e:02} with
$$W(x)\geq Cc_2^{q-(p-1)}|x|^{{\tilde\gamma_-}(q-p+1)+p-\sigma}$$
and $\tilde\gamma_-(q-p+1)+p-\sigma>0$.
Then the assertion follows by Corollary \ref{c:NonExist}.
\end{proof}

\paragraph{Case $\mu=C_H$.}
In this case the proof of the nonexistence can be performed in one step
for both subcritical and critical cases.

\begin{proposition}
Let  $\sigma\leq\gamma_\ast(q-p+1)+p$.
Then \eqref{e:MAIN} has no positive super-solution in $B^c_1$.
\end{proposition}
\begin{proof}
Let $u>0$ be a super-solution to \eqref{e:MAIN} in $B^c_1$.
Then $u$ is a super-solution to
\begin{eqnarray*}
-\Delta_p u-\frac{C_H}{|x|^p}u^{p-1}= 0\quad \mbox{in}\quad B^c_1.
\end{eqnarray*}
By Theorem \ref{l:WCP}(ii) we conclude that $u\geq c|x|^{\gamma_\ast}$ in $B^c_2$.
So $u$ is a super-solution to
\begin{equation}\label{e:04}
-\Delta_p u -\frac{C_H+W(x)}{|x|^p}u^{p-1}=0\quad \mbox{in}\quad B^c_2,
\end{equation}
where
$$W(x):=C|x|^{p-\sigma}u^{q-(p-1)}\geq Cc^{q-(p-1)}|x|^{\gamma_\ast(q-p+1)+p-\sigma},$$
with $\gamma_\ast(q-p+1)+p-\sigma\geq 0$.
Then the assertion follows by Corollary \ref{c:NonExist}.
\end{proof}

\subsection{A nonlinear lower bound}

We will use the comparison principle (Theorem \ref{l:WCPimpr} 
in order to establish the following lower bound on positive solutions
to nonlinear equation \eqref{e:MAIN} in the sub-homogeneous case $q<p-1$.

\begin{lemma}\label{l:EstSubLin}
Let $q<p-1$.
Let $u>0$ be a solution to \eqref{e:MAIN} in $B^c_1$.
Then there exists $c>0$ such that
\begin{equation}\label{e:EstSubLin}
u\geq c|x|^{\frac{\sigma-p}{q-(p-1)}}\quad\mbox{in }\:B^c_{2}.
\end{equation}
\end{lemma}
\begin{proof}
Let $u>0$ be a solution to \eqref{e:MAIN} in $B^c_1$.
Let $x=Ry$ with $y\in A_{2,R}$ and $R\ge 1$. Set
$$v_R(y):=R^{-\frac{\sigma-p}{q-(p-1)}}u(Ry).$$
Then $v_R(y)$ satisfies
\begin{equation}\label{ppe:MAIN}
-\Delta_p v_R
-\frac{\mu}{|y|^p}{v_R}^{p-1}=\frac{C}{|y|^\sigma}{v_R}^q\quad\mbox{in}\quad A_{2,4}.
\end{equation}
Let $\lambda_1>0$ be the principal eigenvalue and $\phi_1>0$ be the principal eigenfunction to
$$-\Delta_p \phi
-\frac{\mu}{|y|^p}{\phi}^{p-1}=\lambda\phi,\qquad\phi\in W^{1,p}_0(A_{2,4}).
$$
By a direct computation, $\tau_0\phi_1$ is a sub-solution to \eqref{e:MAIN}
for a sufficiently small $\tau_0>0$.
Therefore, Theorem \ref{l:WCPimpr} implies that
$$v_R\ge\tau_0\phi_1\quad\mbox{in }\:A_{2,4}.$$
So, lower bound \eqref{e:EstSubLin} follows.
\end{proof}

\subsection{Nonexistence: sub-homogeneous case $q<p-1$}
As before, we distinguish the cases $\mu<C_H$ and $\mu=C_H$.

\paragraph{Case $\mu<C_H$.}
First we consider the subcritical case,
when $(q,\sigma)$ is below to the critical line $\Lambda_\ast$.

\begin{proposition}\label{p:1}
Let  $\sigma<\gamma_+(q-p+1)+p$.
Then \eqref{e:MAIN} has no positive super-solution in $B^c_1$.
\end{proposition}
\begin{proof}
Let $u>0$ be a super-solution to \eqref{e:MAIN} in $B^c_1$. According to Lemma \ref{l:solSubSuper},
we may assume that $u$ is a solution to \eqref{e:MAIN} in $B^c_1$.
Then $u$ is a super-solution to the homogeneous equation
\begin{equation}\label{e:-2}-\Delta_p u-\frac{\mu}{|x|^p} u^{p-1}= 0\quad \mbox{in}\quad B^c_1.\end{equation}
By Theorem \ref{l:upper_est<C_H}(i) we conclude that
\begin{equation}\label{e:33'}
\inf_{S_R}u\leq c_1 R^{\gamma_+},\qquad R>2.
\end{equation}
Since $\gamma_+<\frac{\sigma-p}{q-(p-1)}$ this contradicts to lower bound \eqref{e:EstSubLin}.
\end{proof}

Next we prove the nonexistence in the critical case.
When $(q,\sigma)$ belongs to the critical line $\Lambda_\ast$,
\eqref{e:33'} is no longer incompatible with \eqref{e:EstSubLin},
so we need to improve estimate \eqref{e:33'}.

\begin{proposition}\label{p:3}
Let  $\sigma=\gamma_+(q-p+1)+p$.
Then \eqref{e:MAIN} has no positive super-solution in $B^c_1$.
\end{proposition}
\begin{proof}
Let $u>0$ be a super-solution to (\ref{e:MAIN}) in $B^c_1$.
According to Lemma \ref{l:solSubSuper},
we may assume that $u$ is a solution to \eqref{e:MAIN} in $B^c_1$.
Using \eqref{e:EstSubLin} we conclude that $u>0$ is a solution to
\begin{equation}\label{e:-3}-\Delta_p u-\frac{\mu+W(x)}{|x|^p}u^{p-1}=0\quad\mbox{in}\quad B^c_1,\end{equation}
where $W(x):=C|x|^{p-\sigma}u^{q-(p-1)}\in L^\infty(B^c_2)$. Thus
the strong Harnack Inequality \eqref{pe:sHarnack} combined with
upper bound \eqref{e:33'} implies that
$$\sup_{A_{R/2,R}}u\leq C_S\inf_{A_{R/2,R}}u\leq cR^{\gamma_+},\qquad R>4,$$
and hence $W(x)\geq c_1$ in $B^c_4$, for some $c_1>0$.
Therefore $u>0$ is a super-solution to
\begin{equation}\label{e:34''}
-\Delta_p u-\frac{\tilde\mu}{|x|^p}u^{p-1}=0\quad\mbox{in}\quad B^c_4,
\end{equation}
where $\tilde{\mu}=\mu+\epsilon$ with $0<\epsilon<c_1$ small enough.
Without loss of generality we may assume that $\tilde{\mu}<C_H$.
Then by Theorem~\ref{l:WCP}(i) we conclude that
$\inf_{S_R}u\leq cR^{\tilde\gamma_+}$ for all $R>4$,
where $\tilde\gamma_+\in(\gamma_\ast,\gamma_+)$ is the largest  root of the equation \eqref{e:ROOTS-A} with
$\tilde{\mu}$ in place of $\mu$.
This improved estimate contradicts to lower bound \eqref{e:EstSubLin}.
\end{proof}

\paragraph{Case $\mu=C_H$.}
First we prove the nonexistence in the subcritical case,
when $(q,\sigma)$ is below to the critical line $\Lambda_\ast$.

\begin{proposition}\label{p:4}
Let $\sigma<\gamma_\ast(q-p+1)+p$.
Then \eqref{e:MAIN} has no positive super-solution in $B^c_1$.
\end{proposition}
\begin{proof}
We start as in the proof of Proposition \ref{p:1} with $C_H$ in place of $\mu$ in \eqref{e:-2}.
By Theorem \ref{l:upper_est<C_H}(ii) we conclude that
\begin{equation}\label{e:33''}
\inf_{S_R}u\leq c R^{\gamma_\ast}(\log R)^{\beta^*},\quad R>2,
\end{equation}
where $\beta^*=1$ for $p=N$ and $\beta^*=\frac{2}{p}$ for $p\neq N$.
This contradicts to lower bound \eqref{e:EstSubLin}.
\end{proof}

Now we consider the critical case,
i.e.\ when $(q,\sigma)$ belongs to the critical line $\Lambda_\ast$.
We need to distinguish between the cases $q>-1$ and  $q=1$.

\begin{proposition}\label{p:2}
Let $q\in(-1,p-1)$ and $\sigma=\gamma_\ast(q-p+1)+p$.
Then \eqref{e:MAIN} has no positive super-solution in $B^c_1$.
\end{proposition}
\begin{proof}
We start as in  Proposition \ref{p:3} with $C_H$ in place of $\mu$ in \eqref{e:-3}.
The strong Harnack Inequality \eqref{pe:sHarnack} and upper bound
\eqref{e:33''} imply that
\begin{eqnarray}\label{e:46}
\sup_{A_{R/2,R}}u\leq C_S\inf_{A_{R/2,R}}u\leq cR^{\gamma_\ast}(\log R)^{\beta^*},
\quad R>4.
\end{eqnarray}
We conclude that
$$W(x)\geq\epsilon(\log |x|)^{\beta^*(q-p+1)}\quad\mbox{in }\:B^c_{4},$$
for some $\epsilon>0$.
Hence $u>0$ is a super-solution to the equation
\begin{equation}\label{e:34'}
-\Delta_p u-\frac{C_H}{|x|^p}u^{p-1}-\frac{\epsilon(\log |x|)^t}{|x|^p\log^{m_*}|x|}u^{p-1}=0
\quad\mbox{in}\quad B^c_{4},
\end{equation}
where $t:=\beta^*(q-p+1)+m>0$.
So, the assertion follows by Corollary \ref{c:NonExist}.
\end{proof}

In the 'double critical' case $q=-1$ equation \eqref{e:34'} does not directly lead to the nonexistence, because $t=0$.
So we need to improve estimate \eqref{e:46}.

\begin{proposition}
Let $q=-1$  and $\sigma=\gamma_\ast(q-p+1)+p$.
Then \eqref{e:MAIN} has no positive super-solution in $B^c_1$.
\end{proposition}
\begin{proof}
Arguing as in the proof of Proposition \ref{p:2},
we conclude that $u>0$ is a super-solution to the equation \eqref{e:34'} with t=0.
We may assume that $\epsilon_1<C_\ast$.
Then using Lemma~\ref{l:upper_est<C_H}(iii) and applying the strong Harnack inequality
to equation \eqref{e:34'}, we conclude that
\begin{eqnarray*}
\sup_{A_{R/2,R}}u\leq C_S\inf_{A_{R/2,R}}u\leq cR^{\gamma_\ast}(\log R)^{\beta},
\quad  R>2\rho,
\end{eqnarray*}
where $\beta\in(\beta_+,\beta_*)$ and $\rho> 4$.
Therefore $u>0$ is a super-solution to the equation
$$-\Delta_p u-\frac{C_H}{|x|^p}u^{p-1}
-\frac{W(x)}{|x|^p\log^{m_\ast}|x|}u^{p-1}=0
\quad\mbox{in}\quad B^c_{2\rho},$$
where
$$W(x):=C|x|^{p-N}\log^{m_\ast}|x|\,u^{-p}\geq c(\log|x|)^{2-\beta p}
\quad\mbox{in }\:B^c_{2\rho}.$$
Hence, the assertion follows by Corollary \ref{c:NonExist}.
\end{proof}

This completes the  description of the nonexistence region $\mathcal N$
and the proof of the nonexistence part of Theorem \ref{t:MAIN}.
Next we show that the established nonexistence results are sharp.

\subsection{Existence}

As soon as the nonexistence region $\mathcal N$ is described,
the construction of explicit super-solutions in its complement is straightforward.

\paragraph{Case $\mu<C_H$.}
Let $(q,\sigma)\in \R^2\setminus\mathcal N$.
Choose $\gamma\in(\gamma_-,\gamma_+)$ such that
$$\left\{
\begin{array}{ccl}
\gamma_-<\gamma<\frac{\sigma-p}{q-p+1}&\text{if}&q>p-1,\smallskip\\
\frac{\sigma-p}{q-p+1}<\gamma<\gamma_+&\text{if}&q<p-1,\smallskip\\
\gamma_-<\gamma<\gamma_+&\text{if}&q=p-1.\\
\end{array}
\right.$$
Then one can verify directly that the functions
$u=\tau r^\gamma$ are super-solutions to \eqref{e:MAIN} in $B^c_\rho$
for an appropriate choice of $\tau>0$ and $\rho\ge 1$.

\paragraph{Case $\mu=C_H$.}
Let $(q,\sigma)\in \R^2\setminus\mathcal N$.
For $p=N$, choose $\beta\in(0,1)$ such that
$$\left\{
\begin{array}{ccl}
0<\beta<1&\text{if}&\sigma>N,\quad q\in\R,\smallskip\\
\frac{N}{N-1-q}<\beta<1&\text{if}&\sigma=N,\quad q<-1.\\
\end{array}
\right.$$
Then one can verify directly that the functions
$u=\tau \log^\beta r$
are super-solutions to \eqref{e:MAIN} in $B^c_\rho$
for an appropriate choice of $\tau>0$ and $\rho>1$.

For $p\neq N$, choose $\beta\in(0,2/p)$ such that
$$\left\{
\begin{array}{ccl}
0<\beta<\frac{2}{p}&\text{if}&\sigma>\Lambda_\ast(q),\quad q\in\R,\smallskip\\
-\frac{2}{q-p+1}<\beta<\frac{2}{p}&\text{if}&\sigma=\gamma_\ast(q-p+1)+p,\quad q<-1.\\
\end{array}
\right.$$
Then \eqref{e:A} implies that the function $u=\tau r^{\gamma_\ast}(\log r)^\beta$ satisfies
\begin{equation}\label{e:42}
-\Delta_p u-\frac{C_H}{|x|^p}u^{p-1}
\ge\frac{\epsilon}{|x|^{p}\log^2|x|}u^{p-1}\ge \frac{C}{|x|^\sigma}u^q
\quad\mbox{in }\:B^c_{\rho},
\end{equation}
where $\epsilon=\beta(p-1)(2-\beta p)/2\in(0,C_*)$ and $\tau>0$, $\rho>1$
are  chosen appropriately.

This completes the proof of Theorem \ref{t:MAIN}.

\newpage
\appendix

\section{Appendix: Sample sub- and super-solutions}\label{A1}

\begin{small}
In this section we construct explicit super- and sub-solutions to the homogeneous  equation of the form
\begin{equation}\label{e:44}-\Delta_p
u-\frac{\mu}{|x|^p}u^{p-1}-\frac{\epsilon}{|x|^p\log^{m_*}|x|}u^{p-1}=0\quad
\mbox{in}\quad B^c_\rho,
\end{equation}
where $\rho\geq 2$.
In what follows we assume that $\mu\leq C_H$ and $\epsilon\in[0,C_*)$,
where  $C_H$,  $C_*$ and $m_*$ are defined in \eqref{E:HardyNotations}.
When $u$ is radially symmetric we loosely write $u(|x|)=u(r)$ instead of $u(x)$. In this case
in the
polar coordinates $(r,\omega)$ on $\R^N$ equation \eqref{e:44}
transforms into the ordinary differential equation
\begin{equation}\label{e:rad}
-r^{1-N}(r^{N-1}|u_r|^{p-2}u_r)_r-\frac{\mu}{r^p}u^{p-1}-\frac{\epsilon}{r^p\log^{m_*} r}u^{p-1}=0
\qquad(r>\rho).
\end{equation}
Let $\mu\leq C_H$.
Set $\gamma_\ast:=\frac{p-N}{p}$.
By $\gamma_-\le\gamma_+$ we denote the real roots of the equation
\begin{equation}
\label{e:ROOTS-A}
-\gamma|\gamma|^{p-2}\left(\gamma(p-1)+N-p\right)=\mu.
\end{equation}
If $\mu<C_H$ then $\gamma_-<\gamma_\ast<\gamma_+$.
If $\mu=C_H$ then $\gamma_\pm=\gamma_\ast$.
It is straightforward to see that if $\mu\le C_H$ and $\epsilon=0$ then
{\sl the function $u=r^{\gamma}$ is a sub-solution to equation \eqref{e:rad}
if $\gamma\in(-\infty,\gamma_-]\cup[\gamma_+,+\infty)$ and a
super-solution if $\gamma\in[\gamma_-,\gamma_+]$}.

Let $p=N$ and $\epsilon\in[0,C_*]$.
Then $\beta_-\le\beta_+$ denote the real roots of the equation
\begin{equation}\label{e:80 N}
\beta^{N-1}(1-\beta)(N-1)=\epsilon.
\end{equation}
Notice that $0\le\beta_-\le\frac{N-1}{N}\le\beta_+\le 1$.
It is simple to verify that {\sl the function $u:=\log^{\beta}r$
is a sub-solution to \eqref{e:rad} if
$\beta\in(-\infty,\beta_-]\cup[\beta_+,+\infty)$ and a
super-solution if $\beta\in[\beta_-,\beta_+]$}.

When $p\neq N$ , $\mu=C_H$ and $\epsilon\in[0,C_\ast]$ the situation becomes more delicate.
We denote by $\beta_-\le\beta_+$ the real roots of the equation
\begin{equation}\label{e:80><N}
\frac{1}{2}|\gamma^\ast|^{p-2}(p-1)(2-\beta p)\beta=\epsilon,
\end{equation}
If $\epsilon<C_*$ then  $0\le\beta_-<\frac{1}{p}<\beta_+\le\frac{2}{p}$.
If $\epsilon=C_*$ then  $\beta_-=\beta_+=\frac{1}{p}$.

\begin{lemma}\label{p:SubSuper'}
Let $p\neq N$, $\mu=C_H$ and $\eps\in[0,C_\ast]$.
Let $u_{\beta,\tau}(r):=r^{\gamma_\ast}(\log r)^\beta(\log\log r)^\tau$,
where $\beta\ge 0$ and $\tau\in\R$.
The following assertions are valid.
\begin{enumerate}
\item[{$(i)$}]
Let $\epsilon\in[0,C_*)$.
Then there exists $\rho=\rho(p,N,\beta,\tau)>1$ such that
\begin{enumerate}
\item[{$(a)$}]
$u_{\beta,\tau}$ is a super-solution to \eqref{e:rad}
if $\beta\in(\beta_-,\beta_+)$ and a sub-solution if
$\beta<\beta_-$ or $\beta>\beta_+$;
\item[{$(b)$}]
$u_{\beta_-,\tau}$ is a super-solution to \eqref{e:rad}
if $\tau>0$ and a sub-solution if $\tau<0$;
\item[{$(c)$}]
$u_{\beta_+,\tau}$ is a super-solution to \eqref{e:rad}
if $\tau<0$ and a sub-solution if $\tau>0$;
\item[{$(d)$}]
$u_{\beta_\pm,0}$ is a super-solution to \eqref{e:rad} if $p\in(1,2]\cup(N,+\infty)$
and a sub-solutions if $p\in[2,N)$.
\end{enumerate}
\item[{$(ii)$}]
Let $\epsilon=C_*$.
Then  there exists $\rho=\rho(p,N,\beta,\tau)>1$ such that
\begin{enumerate}
\item[{$(a)$}]
$u_{\beta,\tau}$ is a sub-solution to \eqref{e:rad} if $\beta\neq 1/p$;
\item[{$(b)$}]
$u_{1/p,\tau}$ is a super-solution to \eqref{e:rad} if $\tau\in(0,\frac{2}{p})$
and a sub-solution if $\tau<0$ or $\tau>2/p$;
\item[{$(c)$}]
$u_{1/p,0}$ is a super-solution to \eqref{e:rad} if $p\in (1,2]\cup(N,+\infty)$
and a sub-solution if $p\in[2,N)$.
\end{enumerate}
\end{enumerate}
\end{lemma}

\begin{proof}Observe  that for every $\beta,\tau\in\R$ there exists $\rho>2$
such that $u_r$ does not change sign on $(\rho,\infty)$.
Then a direct computation (similar to  \cite[Lemmas 2.1, 2.2]{Sandeep})
verifies that
\begin{eqnarray}
-\Delta_p\,u_{\beta,\tau}&=&
|\gamma_\ast|^{p-2}\left(|\gamma_\ast|^2+
\frac{\beta(p-1)(2-\beta p)}{2}\frac{1}{\log^2 r}+
\tau(p-1)(1-\beta p)\frac{1}{\log^2 r\log\log r}\right.\nonumber\\\label{e:A}
&+&
\left.\frac{\tau(p-1)(2-\tau p)}{2}\frac{1}{\log^2 r\,\log\log^2 r}
+\frac{p(p-1)(p-2)}{3(N-p)}\beta^2(\beta p-3)\frac{1}{\log^3 r}
+R(r)\right)\frac{u^{p-1}}{r^p},
\end{eqnarray}
where
$$R(r)=O\left(\frac{1}{\log^3 r\log\log r}+\frac{1}{\log^4 r}\right)
\quad\mbox{as }\:r\to\infty.$$
The rest of the proof is straightforward.
\end{proof}

\begin{remark}
Table \ref{table} summarizes some values of the parameters $\beta,\tau\in\R$
which make the function $u_{\beta,\tau}=r^{\gamma_\ast}(\log r)^\beta(\log\log r)^\tau$
a sub- or super-solution to \eqref{e:rad} with $\mu=C_H$ and $\epsilon\in[0,C_*]$,
for a sufficiently large radius $\rho>1$. Observe that the radius $\rho>1$
depends on the data and, in general, can not be determined explicitly.
Similar calculations with $\tau=0$ were provided in \cite{Sandeep,SandeepSreenadh}
for interior domains.
\end{remark}

\begin{table}[t]
\begin{center}
\begin{footnotesize}
\begin{tabular}{|l|c|c|c|}
\hline
& {\sl Sub-solution} & {\sl Super-solution} & {\sl Sub-solution}\\
\hline
$p=N$, $\eps=0$ &
$\beta\le \beta_-\,$, $\tau=0$ &
$\beta\in[\beta_-,\beta_+]$, $\tau=0$ &
$\beta\ge\beta_+\,$, $\tau=0$\\
\hline
$p\neq N$, $\eps=0$ &
$\beta\le 0\,$, $\tau=0$ &
$\beta\in[0,2/p]$, $\tau=0$ &
$\beta>2/p\,$, $\tau=0$\\
&
&
&
$\beta=2/p\,$, $\tau>0$\\
\hline
$p\neq N$, $\eps\in(0,C_\ast)$ &
$\beta<\beta_-\,$, $\tau=0$ &
$\beta\in(\beta_-,\beta_+)$, $\tau=0$ &
$\beta>\beta_+\,$, $\tau=0$\\
&
$\beta=\beta_-\,$, $\tau<0$ &
$\beta=\beta_-$, $\tau>0$ or $\beta=\beta_+$, $\tau<0$ &
$\beta=\beta_+\,$, $\tau>0$\\
\hline
$p\neq N$, $\eps=C_\ast$ &
$\beta<1/p\,$, $\tau=0$ &
&
$\beta>1/p\,$, $\tau=0$\\
&
$\beta=1/p\,$, $\tau<0$ &
$\beta=1/p\,$, $\tau\in(0,2/p)$&
$\beta=1/p\,$, $\tau>2/p$\\
\hline
\end{tabular}
\caption{{\small\em Case $\mu=C_H$.
Properties of
$u_{\beta,\tau}=r^{\gamma_\ast}(\log r)^\beta(\log\log r)^\tau$,
for a large $\rho>1$.}}\label{table}
\end{footnotesize}
\end{center}
\end{table}

\section{Appendix: Small sub-solutions }

A {\sl small} (sub)\,solution to equation \eqref{e:44}
is a (sub)\,solution $v>0$ to \eqref{e:44} that satisfies the condition:
\begin{enumerate}
\item[$(\mathcal{S})$]
there exists a sequence $(\theta_n)\in W_c^{1,\infty}(\R^N)$ such that
$\theta_n\to 1$ a.e. in $\R^n$ and
$$\int_{B^c_\rho}\mathcal{R}(\theta_nv,v)\,dx \to 0\quad\mbox{as }\:R\to +\infty,$$
\end{enumerate}
where $\mathcal{R}(w,v)=
|\nabla w|^p-\nabla\left(\frac{w^p}{v^{p-1}}\right)|\nabla v|^{p-2}\nabla v$
is defined as in Proposition \ref{p:Picon}.
In order to apply Theorem \ref{l:WCPimpr} to equation \eqref{e:44}
we need to verify that \eqref{e:44} has small sub-solutions, which is done in the following proposition.

\begin{proposition}\label{p:R(v,v)}
Set $v=r^\gamma(\log r)^\beta (\log\log r)^\tau$.
The following assertions are valid.
\begin{enumerate}
\item[{$(i)$}]
Let $\gamma\leq\gamma_\ast$, $\beta=0$, $\tau=0$.
Then $v$ is a small sub-solution to \eqref{e:44} with $\mu\leq C_H$ and $\epsilon =0$;
\item[{$(ii)$}]
Let $p\neq N$, $\gamma=\gamma_\ast$, $\beta=1/p$, $\tau<0$.
Then $v$ is a small sub-solution to \eqref{e:44} with $\mu= C_H$ and $\epsilon \in(0, C_*)$;
\item[{$(iii)$}]
Let $p=N$, $\gamma=\gamma_\ast$, $\beta=\frac{N-1}{N}$, $\tau<0$.
Then $v$ is a small sub-solution to \eqref{e:44} with $\mu=0$ and $\epsilon\in(0,C_*)$.
\end{enumerate}
\end{proposition}

\begin{proof}
Lemma \ref{p:SubSuper'} in  case (ii) and  direct computations in cases (i), (iii) show that $v$
is a sub-solution to \eqref{e:44} for corresponding $\mu$ and $\epsilon$.
Below we show that
$\int_{B^c_\rho}\mathcal{R}(\theta_R^\alpha v,v)\,dx\rightarrow 0$ as $R\to+\infty,$
where $\theta_R\in C^{1,1}(0,\infty)$ is defined by
$$\theta_R(r):=\left\{
\begin{array}{cl}
1, & 0\leq r\leq R,\smallskip\\
\frac{\log\frac{R^2}{r}}{\log R},& R\leq r\leq R^2,\smallskip\\
0, & r\geq R,
\end{array}
\right.$$
and $\alpha\geq 1$ will be chosen later.
By Proposition \ref{p:Picon}, for $R>\rho$ we have
\begin{eqnarray*}\int_{B_\rho^c}\mathcal{R}(\theta_R^\alpha v,v)\,dx &=&
\int_{A_{\rho,R}}\mathcal{R}(\theta_R^\alpha v,v)\,dx
+ \int_{A_{R,R^2}}\mathcal{R}(\theta_R^\alpha v,v)\,dx\\&=&
\int_{A_{\rho,R}}\mathcal{R}(v,v)\,dx
+ \int_{A_{R,R^2}}\mathcal{R}(\theta_R^\alpha v,v)\,dx
=c_N\int_{R}^{R^2}\mathcal{R}(\theta_R^\alpha v,v)r^{N-1}dr
\end{eqnarray*}
Below we estimate the latter integral.

\noindent
$(i)$ Using the inequalities
(see, e.g., \cite[Lemma 7.4]{Shafrir} )
\begin{eqnarray}
\mathcal R(\theta_R v,v)&\le& c_1|\theta_R v_r'|^{p-2}|v(\theta_R)_r'|^2+c_2|v\,(\theta_R)_r'|^p,
\qquad(p>2),\label{e-sub>2}\\
\mathcal R(\theta_R v ,v)&\le& c_3|v\,(\theta_R)_r'|^p,
\hspace{9em}\qquad (1<p\le 2),\label{e-sub<2}
\end{eqnarray}
we obtain directly that there exists $c>0$ such that
\begin{equation}\label{B1-i}
\int_{R}^{R^2}\mathcal{R}(\theta_R v,v) r^{N-1}dr\le c\frac{R^{\gamma p+N-p}}{(\log R)^p}.
\end{equation}

\noindent
$(ii)$
Set $Q(r):=-\gamma_\ast\log r\log\log r-\beta\log\log r-\tau.$
Then direct
computations give
\begin{eqnarray*}
\mathcal{R}(\theta^\alpha_R v,v)&=&\frac{(\log r)^{(\beta-1)p}
(\log\log r)^{(\tau-1)p} }{r^N (\log R)^{\alpha p}}\,
\left(\log\frac{R^2}{r}\right)^{\alpha p-p}\,Q^p(r)\times\\
&\times&\left\{\left|\, \log\frac{R^2}{r}+\alpha\,\frac{\log r\log\log r}{Q(r)}
\right|^p-\left(\log\frac{R^2}{r}\right)^{p-1}\left(\log\frac{R^2}{r}+\alpha
p\frac{\log r\log\log r}{Q(r)}\right) \right\}.
\end{eqnarray*}

\noindent
Let $p\geq 2$. Choose $\alpha=1$.  We use  the inequality
(see, e.g., \cite[Lemma 7.4]{Shafrir})
\begin{equation}\label{e:z>2<2}
\quad|z_1+z_2|^p-|z_1|^p-p|z_1|^{p-2}z_1z_2\leq
\frac{p(p-1)}{2}\left(|z_1|+|z_2|\right)^{p-2}|z_2|^2\qquad(z_1,z_2\in\R)
\end{equation}
with
$$z_1=\log\frac{R^2}{r},\qquad z_2=\frac{\log r\log\log r}{Q(r)}$$
to obtain that
\begin{eqnarray*}
\mathcal{R}(\theta_R v,v)
&\leq&
c_1\,\frac{(\log r)^{(\beta-1)p+2}(\log\log r)^{(\tau-1)p+2}}{r^N \log^p R}
\,Q^{p-2}(r)\left|\,\log\frac{R^2}{r}+\frac{\log r\log\log r}{Q(r)}\right|^{p-2}\\
&=&
c_1\,\frac{(\log r)^{(\beta-1)p+2}(\log\log r)^{(\tau-1)p+2}}{r^N\log^p R}
\,\left|\,Q(r)\,\log\frac{R^2}{r}+\log r\log\log r\right|^{p-2}.
\end{eqnarray*}
Thus we arrive at
\begin{equation}\label{B1-ii}
\int_{A_{R,R^2}}\mathcal{R}(\theta_R v,v)\,dx \leq
\frac{c_1}{\log^2 R}
\int_{R}^{R^2}\frac{(\log r)^{\beta p}(\log\log r)^{\tau p}}{r}\,dr
\le c(\log R)^{\beta p-1}(\log\log R)^{\tau p}.
\end{equation}
\smallskip

\noindent
If $1<p<2$, choose $\alpha> \frac{2}{p}$. Observe that the  Taylor expansion 
applied to the function $f(t)=|z_1+tz_2|^p$ with
$0<t<1$, $z_1,z_2\in\R$, $z_1\neq 0$ and $z_1 z_2\ge 0$
leads to
$$|z_1+z_2|^p-|z_1|^p-p|z_1|^{p-2}z_1z_2=
\frac{p(p-1)}{2}\,|z_1+t_0 z_2|^{p-2}|z_2|^2\leq
\frac{p(p-1)}{2}\,|z_1|^{p-2}|z_2|^2,\quad t_0\in (0,1).$$ Using the above inequality
with
$$z_1=\log\frac{R^2}{r},\qquad z_2=\alpha\,\frac{\log r\log\log r}{Q(r)},$$
we obtain
\begin{eqnarray*}
\mathcal{R}(\theta_R v,v)&\leq&
c_1\,\frac{(\log r)^{(\beta-1)p+2}(\log\log r)^{(\tau-1)p+2}}{r^N (\log R)^{\alpha p}}\,
\left(\log\frac{R^2}{r}\right)^{\alpha p-2}\,Q^{p-2}(r).
\end{eqnarray*}
Since $\alpha p-2>0$ we conclude that
\begin{equation}\label{B1-ii-2}
\int_{A_{R,R^2}}\mathcal{R}(\theta_R v,v)\,dx \leq
\frac{c_1}{\log^2 R}\int_{R}^{R^2}\frac{(\log r)^{\beta p}(\log\log r)^{\tau p}}{r}\,dr
\le c(\log R)^{\beta p-1}(\log\log R)^{\tau p}.
\end{equation}

\noindent
$(iii)$ An easy computations shows that
$$\mathcal{R}(\theta_R v,v)\leq \frac{c}{r^N\log^N R}(\log r)^{\beta N}(\log\log r)^{\tau N}.$$
Therefore
\begin{eqnarray}\label{B1-iii}
\int_{A_{R,R^2}}\mathcal{R}(\theta_R v,v)\,dx &\leq&
\frac{c}{\log^N R}\int _{R}^{R^2}
\frac{(\log r)^{\beta N}(\log\log r)^{\tau N}}{r}\,dr\nonumber\\
&\le& c(\log R)^{\beta N-N+1}(\log\log r)^{\tau N}.
\end{eqnarray}
This completes the proof.
\end{proof}

\section{Appendix: Large sub-solutions}\label{A:B}

A {\sl large} (sub)\,solution to equation \eqref{e:44}
is a positive (sub)\,solution of the problem
\begin{equation}\label{e:66}
-\Delta_p u-\frac{\mu}{|x|^p}u^{p-1}-\frac{\epsilon}{|x|^p\log^{m_\ast}|x|}u^{p-1}(\,\leq\,)=0\quad\mbox{in }\:B^c_R,
\qquad u=0\quad\mbox{on }\:S_R,
\end{equation}
with a sufficiently large $R>1$.
Below we establish the existence and asymptotic behavior  of large sub-solutions.

\begin{theorem} \label{t:B0}
Let $\mu\le 0$ and $\epsilon=0$. The following assertions are valid.
\begin{enumerate}
\item[$(i)$] if $p\neq N$ or $\mu<0$ then $u=|x|^{\gamma_+}-R^{\gamma_+}$ is a positive sub-solution to  \eqref{e:66}.
\item[$(ii)$]
if $\mu=0$ and $p=N$ then
$u=\log|x|-\log R$ is a positive sub-solution to  \eqref{e:66}.
\end{enumerate}
\end{theorem}

\begin{proof}
Note that if $\mu\le 0$ then $0\in[\gamma_-,\gamma_+]$. Hence
positive constants are super-solutions to \eqref{e:66}.
Then a direct computation verifies that
$u=r^{\gamma_+}- R^{\gamma_+}$ or $u=\log|x|-\log R$   are sub-solutions to \eqref{e:66}.
\end{proof}

\begin{theorem} \label{t:B1}
The following assertions are valid.
\begin{enumerate}
\item[$(i)$]
Let $p\neq N$, $\mu\in (0,C_H)$ and $\epsilon=0$.
Then \eqref{e:66} admits a solution $u>0$ such that
$$u=c|x|^{\gamma_+(1+o(1))}\quad\mbox{as}\quad r\to+\infty.$$
\item[$(ii)$]
Let $p\neq N$, $\mu=C_H$ and $\epsilon=0$.
Then \eqref{e:66} admits a solution $u>0$ such that
$$u=c|x|^{\gamma_\ast}\left(\log|x|\right)^{\frac{2}{p}\,(1+o(1))}\quad\mbox{as}\quad r\to+\infty.$$
\item[$(iii)$]
Let $p=N$, $\mu=C_H$ and $\epsilon\in(0,C_\ast)$.
Then \eqref{e:66} admits a solution $u>0$ such that
$$u=c\left(\log|x|\right)^{\beta_+(1+o(1))}\quad\mbox{as}\quad r\rightarrow+\infty.$$
\item[$(iv)$]
Let $p\neq N$, $\mu=C_H$ and $\epsilon\in(0,C_*)$.
Then \eqref{e:66} admits a solution $u>0$ such that for every
$\delta\in(0,\min\{\beta_+-\frac{1}{p},\frac{2}{p}-\beta_+\})$ there exists $c_\delta>0$ and $R_\delta>e$ and $u$
satisfies
$$c^{-1}_\delta|x|^{\gamma_\ast}\left(\log|x|\right)^{\beta_+-\delta}\leq
u\leq c_\delta|x|^{\gamma_\ast}\left(\log|x|\right)^{\beta_++\delta}
\quad\mbox{in}\quad (R_\delta,+\infty).$$
\end{enumerate}
\end{theorem}

Our proof of Theorem \ref{t:B1} employs the generalized Pr\"ufer Transformation.
The classical Pr\"ufer transformation is a well-known  tool in the theory
of linear second-order elliptic equations, cf. \cite[Chapter 8]{Cod}.
Its generalization to the context of $p$-Laplace equations was recently introduced
by Reichel and Walter \cite{ReichelWalter}, see also \cite{Binding,Brown-Reichel}.
For the readers' convenience we collect below required facts for the generalized sine functions
and Pr\"ufer transformation.

\subsection{Generalized sine function}

The generalized sine function $S_p(\psi)$ ($p>1$) was introduced in \cite{Lindquist-Prufer}
as the solution to the problem
\begin{equation}\label{e:B1}
|w'|^p+\frac{|w|^p}{p-1}=1,\quad w(0)=0,\quad w'(0)=1.
\end{equation}
Equation \eqref{e:B1} arises as a first integral of
$(w'|w'|^{p-2})'+{w|w|^{p-2}}=0$.
The solution of \eqref{e:B1} defines the function $S_p(\psi)=\sin_p(\psi)$ as long as it is increasing,
that is, for $\psi\in [0,\pi_p/2]$, where
\begin{equation}\label{e:B2}
\frac{\pi_p}{2}=\int_0^{(p-1)^{1/p}}\frac{dt}{1-t^p/(p-1)^{1/p}}\,=\frac{(p-1)^{1/p}}{p\sin(\pi/p)}\,\pi.
\end{equation}
Since $S'_p(\pi_p/2)=0$, we  define $S_p$ on the interval
$[\pi_p/2,\pi_p]$ by $S_p(\psi)=S_p(\pi_p-\psi)$, and for
$\psi\in(\pi_p, 2\pi_p]$ we put $S_p(\psi)=-S_p(2\pi_p-\psi)$ and
extend $S_p$ as a $2\pi_p$ -- periodic function on $\R$.
The following properties of $S_p$ will be used frequently (see \cite{Lindquist-Prufer}).

\begin{lemma}\label{l:sin}
The generalized sine function $S_p$ satisfies the  following properties.
\begin{enumerate}
\item[$(i)$]
$S_p$ satisfies \eqref{e:B1} on $\R$; $S_p\in C^1(\R)$ and $\|S_p\|_\infty=(p-1)^{1/p}$;
\item[$(ii)$]
$S'_p|S'_p|^{p-2}\in C^1(\R)$, $\,\|S'_p\|_\infty=1$ and
$\|(S'_p|S'_p|^{p-2})'\|_\infty=(p-1)^{(p-1)/p}$;
\item[$(iii)$]
if $p\leq 2$ then $S'_p\in C^1(\R)$,
while if $p\geq 2$ then $S'_p\in C^{1,1/(p-1)}(\R)$;
\item[$(iv)$]
$(p-1){S_p}''(\psi)=-S_p^{p-1}|S'_p|^{2-p}$ , $\psi\in(0,\pi_p)$, $\psi\neq \pi_p/2$.
\end{enumerate}
\end{lemma}

Clearly, $S_2(\psi)=\sin(\psi)$ and $\pi_2=\pi$.
Notice also that $S_p(t)\rightarrow 1-|t-1|$ as $p\to\infty$, and $S_p(t)\to 0$ as $p\rightarrow 1$.
The generalized sine function was discussed in great detail by Lindquist in \cite{Lindquist-Prufer}.

\subsection{Generalized Pr\"ufer transformation}
In order to construct a positive solution of \eqref{e:66}
it is sufficient to solve the initial value problem
\begin{equation}\label{e:B3}
-r^{1-N}(r^{N-1}|u_r|^{p-2}u_r)_r-V(r)u^{p-1}=0\quad\mbox{in }\:(R,+\infty),
\qquad u(R)=0,\quad u^\prime(R)>0,
\end{equation}
where we set
$$V(r):=\frac{\mu}{r^p}+\frac{\epsilon}{r^p\log^{m_*}r}.$$
Following \cite{Brown-Reichel},
we use the generalized sine function to transform \eqref{e:B3} into phase space
via the generalized polar coordinates $(\rho,\psi)$ defined by
\begin{equation}\label{e:B4}
\left\{
\begin{array}{lcl}
r^{N-1}u'|u'|^{p-2}&=&\rho(r)S'_p(\psi(r))|S'_p(\psi(r))|^{p-2},\\
Q(r)^{(p-1)/p}u^{p-1}&=&\rho(r)S^{p-1}_p(\psi(r)),
\end{array}
\right.
\end{equation}
where the function $0<Q\in C^1(R,+\infty)$ will bee chosen later.
A calculation similar to \cite[Lemma 2]{ReichelWalter}
shows that by means of the generalized polar coordinates \eqref{e:B4}
equation \eqref{e:B3} transforms  into the Cauchy problem
\begin{equation}\label{e:Bpsi}
\left\{
\begin{array}{lclclcl}
\psi'&=&V_1|S'_p(\psi)|^p+V_2\frac{S^p_p(\psi)}{p-1}
+\frac{1}{p}\frac{Q'}{Q}\,S_p(\psi)S'_p(\psi)|S'_p(\psi)|^{p-2},&&\psi(R)&=&0,\smallskip\\
\rho'&=& \rho\left\{\left(V_1-V_2\right)
S_p^{p-1}(\psi)S'_p(\psi)+\frac{1}{p}\frac{Q'}{Q}S_p^p(\psi)\right\},&&\rho(R)&>&0,
\end{array}
\right.
\end{equation}
in $(R,+\infty)$,
where $V_1$ and $V_2$ are defined by
$$V_1(r):=r^{\frac{1-N}{p-1}}Q^{\frac{1}{p}}(r),\qquad
V_2(r):=r^{N-1}V(r)Q^{\frac{1-p}{p}}(r).$$
Notice also that by means of \eqref{e:B4} a pair of $C^2$--functions $(\rho,\psi)$ satisfying
 \eqref{e:Bpsi} transforms
into a positive solution $u$ to \eqref{e:B3}. 

The main feature of system \eqref{e:Bpsi} is the fact that its first equation is independent of $\rho$.
Notice also that the second equation is linear in $\rho$ and is completely integrable
provided the solution $\psi$ of the first equation is given.

For the choice of $Q(r)$ we distinguish between the cases $V(r)>0$ and $V(r)<0$.
If $V(r)>0$ then we set
\begin{equation}\label{e:Q}
Q(r)=V(r)r^{\frac {p(N-1)}{p-1}}.
\end{equation}
Then $V_1=V_2=V^{1/p}$ and using Lemma \ref{l:sin}
we rewrite \eqref{e:Bpsi} in the form
\begin{equation}\label{e:BpsiV>0}
\left\{
\begin{array}{lclclcl}
\psi'&=&
V^{1/p}+\left(\frac{1}{p}\frac{V'}{V}+\frac{N-1}{p-1}\,\frac{1}{r}\right)S_p(\psi)S'_p(\psi)|S'_p(\psi)|^{p-2},
&&\psi(R)&=&0,\smallskip\\
\rho'&=&
\rho\left(\frac{1}{p}\frac{V'}{V}+\frac{N-1}{p-1}\,\frac{1}{r}\right)S_p^p(\psi),
&&\rho(R)&>&0,
\end{array}
\right.
\end{equation}
in $(R,+\infty)$. In the case $V(r)<0$ one can choose
$Q(r)=-V(r)r^{\frac {p(N-1)}{p-1}}$, however we are not interested in this case below. 

\smallskip
The main tools of our analysis of \eqref{e:BpsiV>0} will be a simple comparison principle between
sub- and super-solutions and a stabilization argument for a time--dependent one-dimensional ODEs.
The comparison principle below can be found in \cite{ReichelWalter}.
\begin{lemma}[{\sf Comparison principle}]\label{l:compar}
Let $f:(R,\infty)\times\R\to\R$ be locally Lipschitz--continuous 
in $(R,\infty)\times\R$.
Let $\phi,\varphi$ be $C^1$--functions on $(R,\infty)$, continuous in $[R,\infty)$, and such that
$$\phi^\prime(r)\le f(r,\phi),\qquad\varphi^\prime(r)\ge f(r,\varphi),
\qquad
\phi(R)\le\varphi(R).$$
Then $\phi(r)\le\varphi(r)$ in $[R,\infty)$.
\end{lemma}

\begin{lemma}[{\sf Stabilization principle}]\label{l:stab}
Let $f:(R,\infty)\times\R\to\R$ be locally  Lipschitz--continuous in $(R,\infty)\times\R$,
and $\lim_{r\to\infty}f(r,\xi)=f_\ast(\xi)$, uniformly on compact subsets of $\R$.
Let $0<\eta\in C^1(R,\infty)$ and $\int_R^{\infty}\eta^{-1}(r)dr=\infty$.
Let $\psi$ be a $C^1$--function on $(R,\infty)$ such that
$$\psi^\prime=\frac{f(r,\psi)}{\eta(r)}\qquad(r>R).$$
Assume that $f(r,\psi(r))>0$ for all $r>R$ and $\psi$ is bounded above.
Then $f_\ast(\psi_\ast)=0$, where $\psi_\ast=\lim_{r\to\infty}\psi(r)$.
\end{lemma}

\begin{proof}
Observe that $\psi(r)$ is monotone increasing and uniformly bounded, so the limit
$\psi_\ast$ exists.  Assume for a contradiction that $f_*(\psi_*)>0$. Then  there exist $\delta>0$ and $R_1>R$ such that
$f(r,\psi(r))>\delta$ for all $r>R_1+1$. Then
$$\psi(r)=\psi(R_1)+\int_{R_1}^r \frac{f(s,\psi(s))}{\eta(s)}ds\ge c_1+\int_{R_1+1}^r\frac{\delta}{\eta(s)}ds
\to\infty\quad\mbox{as }\:r\to+\infty,$$
which contradicts to the boundedness of $\psi$.
Thus the assertion follows.
\end{proof}

\subsection{Proof of Theorem \ref{t:B1}}

Below we establish the existence and asymptotic behavior of
a solution $(\psi,\rho)$ to system \eqref{e:BpsiV>0}.
Then the existence and asymptotic of a positive solution to \eqref{e:B3}
can be computed directly from the asymptotic of $\psi$ and $\rho$ via \eqref{e:B4} and \eqref{e:Q}.

\paragraph{(i) Case $\mu\in(0,C_H)$, $\epsilon =0$, $p\neq N$.}
We consider in detail only the case $p>N$, the case $p<N$ being similar.

System \eqref{e:BpsiV>0} can be written in the form
\begin{equation}\label{e:Bpsi_mu>0}
\psi'=\frac{{F}(\psi)}{r},\qquad
\frac{\rho'}{\rho}=\frac{G(\psi)}{r}\quad\mbox{in}\quad(R,+\infty),
\end{equation}
where
\begin{equation*}
F(\psi):=\mu^{1/p}+\frac{N\!-p}{p-1}\,S_p(\psi)S'_p(\psi)|S'_p(\psi)|^{p-2},\qquad
G(\psi):=\frac{N\!-p}{p-1}\,S_p^p(\psi).
\end{equation*}
%
Notice that $0<\gamma_-<\gamma_+$.
An elementary calculation involving \eqref{e:B1} shows that $F(\psi)=0$ if and only if
$\psi$ satisfies
\begin{equation}\label{e:psi-pm}
S_p(\psi)=\left((\gamma_\pm(p-1)+(N-p))\frac{p-1}{N\!-p}\,\right)^{1/p}
\quad\mbox{and}\quad
S_p'(\psi)=\left(\gamma_\pm\frac{p-1}{p\!-N}\right)^{1/p}.
\end{equation}
Then it follows from the definition and properties of $S_p(\psi)$ that the solutions $\psi_\pm\in(0,\pi_p)$ of
$F(\psi)=0$ are uniquely (modulo $2\pi_p$) determined  by $\gamma_\pm$ via \eqref{e:psi-pm}. One can also see
that
\begin{eqnarray*}
0<\psi_+<\psi_-<\frac{\pi_p}{2}.
\end{eqnarray*}
Moreover, $F(\psi)$ is strictly positive for $\psi\in(0,\psi_+)$.

Let $\psi(r)$ be the solution to the  problem 
\begin{equation}\label{e-psi-mu>0}
\psi'=\frac{{F}(\psi)}{r}\quad\mbox{in}\quad(R,+\infty),\qquad\psi(R)=0,
\end{equation}
for some $R>1$.
Observe that the right hand side of \eqref{e-psi-mu>0} is bounded and smooth
for all $(r,\psi)\in(1,\infty)\times\R$, so $\psi(r)$ exists  for all $r>R$.
Note also that $\psi_+(r)\equiv\psi_+$ is a stationary solution to \eqref{e-psi-mu>0}.
So, $\psi(r)\le\psi_+$ for all $r>R$, by Lemma \ref{l:compar}.
Moreover $\psi(r)$ is monotonically increasing and $F(\psi(r))>0$ for all $r>0$.
Thus, by Lemma \ref{l:stab} we conclude that $\lim_{r\to\infty}\psi(r)=\psi_+$.

\begin{lemma}\label{l:B2} Let $\psi$ be the solution to \eqref{e-psi-mu>0}. Then
$\psi(r)=\psi_++\omega(r)$ where $\omega(r)<0$ in $[R,+\infty)$ and
\begin{eqnarray*}
\omega(r)=c r^{-\left(\gamma_+p+N-p\right)(1+o(1))}\quad
\mbox{as}\quad r\rightarrow +\infty,\end{eqnarray*}
for some $c<0$. \end{lemma}
\begin{proof}
Since 
\begin{eqnarray}\label{e:B19}
F(\psi)=F(\psi_+)+F'(\psi_+)(\psi-\psi_+)+\Theta(\psi-\psi_+)=F'(\psi_+)(\psi-\psi_+)+\Theta(\psi-\psi_+),
\end{eqnarray}
where 
$\Theta(\psi-\psi_+)=O((\psi-\psi_+)^2)$ as $\psi\to\psi_+$,
by Lemma \ref{l:sin} (iv) we obtain that
\begin{equation}\label{e:B22}
F'(\psi)=\frac{N-p}{p-1}\left(|S'_p(\psi)|^p+(p-1)S_p(\psi)|S'_p(\psi)|^{p-2}S''_p(\psi)\right)
=\frac{N-p}{p-1}\left(|S'_p(\psi)|^p-S_p^p(\psi)\right).
\end{equation}
Using \eqref{e:psi-pm} we arrive at $F'(\psi_+)=-(\gamma_+ p+N-p )<0$.
Set $\omega(r):=\psi(r)-\psi_+$. 
Thus $\omega(R)=-\psi_+$ and $\omega$ satisfies
$$\frac{\omega'}{\omega}=\frac{F'(\psi_+)}{r}+\frac{\Theta(\omega)}{r\omega},\qquad r\in(R,+\infty).$$
Therefore we infer that
$$\log\frac{\omega(r)}{\omega(R)}=\log\left(\frac{r}{R}\right)^{-(\gamma_+p+N-p)}+\int_R^r\frac{\Theta(\omega)}{\omega}\,
\frac{ds}{s}.$$
So, the assertion follows by  the L'Hopital's Rule.
\end{proof}

Given the solution $\psi(r)$ to \eqref{e-psi-mu>0}, 
let  $\rho(r)$ be the solution to the problem
\begin{equation}\label{e-rho-mu>0}
\frac{\rho'}{\rho}=\frac{G(\psi)}{r}\quad\mbox{in}\quad(R,+\infty),\qquad\rho(R)=1.
\end{equation}
Observe that the right hand side of \eqref{e-rho-mu>0} is bounded and smooth
for all $(r,\psi)\in(R,\infty)\times\R$, so $\rho(r)$ exists  for all $r>R$.

\begin{lemma}\label{p:mu>0} Let $\rho$  be the solution to \eqref{e-rho-mu>0}. Then
$\rho(r)=cr^{(\gamma_+(p-1)+N-p)(1+o(1))}$ as $r\to+\infty$,
for some $c>0$.
\end{lemma}
\begin{proof} Observe that  $\rho$ satisfies
\begin{equation}\label{e:B-1}\frac{\rho'}{\rho}=\frac{G(\psi_+)}{r}+\frac{\Xi(\omega (r))}{r},\qquad r\in(R,+\infty),
\end{equation}
where  $\Xi(\psi-\psi_+)=o(\psi-\psi_+)$ as $\psi\to\psi_+$ and $\omega(r):=\psi-\psi_+$ is given by Lemma \ref{l:B2}.
Using the definition of $G$, \eqref{e:psi-pm} and \eqref{e:ROOTS-A} we conclude that
$G(\psi_+)=\gamma_+(p-1)+N-p$.
Therefore
$$\log\frac{\rho(r)}{\rho(R)}=\log\left(\frac{r}{R}\right)^{\gamma_+(p-1)+N-p}+\int_R^r\Xi(\omega)\frac{ds}{s}.$$
So, the assertion follows by  the l'Hopital Rule.
\end{proof}

\begin{remark}
The case $\mu\in(0,C_H)$, $\epsilon =0$ and $p<N$ is similar,
the only difference  being  that if $p<N$ then $\gamma_-<\gamma_+<0$
and hence $\pi_p/2<\psi_+<\psi_-<\pi_p$.
\end{remark}

\paragraph{(ii) Case $\mu=C_H$, $\epsilon =0$, $p\neq N$.}
We consider in detail only the case $p>N$,
the case $p<N$ being similar.
System \eqref{e:BpsiV>0} can be written in the form \eqref{e:Bpsi_mu>0},
where
\begin{equation}\label{e:FG-CH}
F(\psi):=|\gamma_\ast|+\frac{N-p}{p-1}\,S_p(\psi)S'_p(\psi)|S'_p(\psi)|^{p-2},\qquad
G(\psi):=\frac{N-p}{p-1}\,S_p^{p}(\psi).
\end{equation}
Notice that $\gamma_*=\frac{p-N}{p}>0$.
A simple analysis shows that $F(\psi)=0$ if and only if
$\psi_\ast=(\pi/4)_p$ modulo $2\pi_p\,$,
where $(\pi/4)_p\in(0,\pi_p/2)$ denotes the unique solution to the equation
\begin{equation}\label{e:B15}
S_p(\psi)=S'_p(\psi)=\left(\frac{p-1}{p}\right)^{1/p}.
\end{equation}
It is clear that $(\pi/4)_2=\pi/4$.
Observe that $F(\psi)$ is nonnegative for all $\psi\in\R$ and
strictly positive for $\psi\in(0,\psi_\ast)$.

Let $\psi(r)$ be the solution to \eqref{e-psi-mu>0},
for some $R>1$.
Clearly $\psi(r)$ exists for all $r>R$.
Note also that $\psi_\ast(r)\equiv\psi_\ast$ is a stationary solution to \eqref{e-psi-mu>0}.
So, $\psi(r)\le\psi_*$ for all $r>R$ by Lemma \ref{l:compar}.
Moreover, by Lemma \ref{l:stab} we conclude that $\lim_{r\to\infty}\psi(r)=\psi_\ast$.

\begin{lemma}\label{l:B2-CH} Let $\psi$ be the solution to \eqref{e-psi-mu>0}. Then
$\psi(r)$ admits a representation $\psi(r)=\psi_\ast+\omega(r)$ where $\omega(r)<0$ in $[R,+\infty)$ and
\begin{eqnarray}\label{e:B29}
\omega(r)=-\frac{2}{p}\,\frac{p-1}{p-N}\,\frac{1+o(1)}{\log r}\quad\mbox{as}\quad r\rightarrow+\infty.
\end{eqnarray}
\end{lemma}
  Proof is similar to  the arguments in the proof of  Lemma \ref{l:B2}. Notice only that
$F'(\psi_*)=F(\psi_*)=0$ and ${F}''(\psi_*)=\frac{(p-N)\,p}{p-1}$, so  use
$F(\psi)=\frac{1}{2}F''(\psi_*)(\psi-\psi_*)^2+o((\psi-\psi_*)^2)$ instead of
 \eqref{e:B19}.

\begin{lemma}\label{p:CH} Let $\rho$  be the solution   to \eqref{e-psi-mu>0}. Then
$\rho(r)=cr^{-\gamma_\ast}(\log r)^{\left(\frac{2(p-1)}{p}\right)(1+o(1))}$
as $r\rightarrow+\infty$, for some $c>0$.
\end{lemma}
Proof is essentially the same as the one for Lemma \ref{p:mu>0}, the only difference being that
instead of \eqref{e:B-1} one uses
$$
\frac{\rho'}{\rho}=\frac{G(\psi_*)}{r}+\frac{G'(\psi_*)\omega(r)}{r}
+\frac{o(\omega(r))}{r},\qquad r\in(R,+\infty),
$$ where $G(\psi_*)=-\gamma_\ast$, $G'(\psi_*)=N-p$.

\begin{remark}
The case $\mu=C_H$, $\epsilon =0$ and $p<N$ is  similar,
the only difference being that 
$\gamma_*<0$ and hence $\psi_\ast=(3\pi/4)_p:=\pi_p-(\pi/4)_p$.
\end{remark}

\paragraph{(iii) Case $\mu=C_H$, $\epsilon \in(0,C_*)$, $p=N$.}
In this case system \eqref{e:BpsiV>0} can be written in the form
\begin{eqnarray*}
\psi'=\frac{{F}(\psi)}{r\log r},\qquad
\frac{\rho'}{\rho}=\frac{G(\psi)}{r\log r}\quad\mbox{in}\quad (R,+\infty),
\end{eqnarray*}
where
\begin{equation*}
F(\psi):=\epsilon^{1/N}
-S_N S_N'|S_N'|^{N-2},\qquad
G(\psi):=-S_N^N.
\end{equation*}
A simple calculation shows that $F(\psi)=0$ if and only if
\begin{equation}\label{e:B30}
S_N(\psi)=(1-\beta_\pm)^{2/N},\qquad S_N'(\psi)=\beta_\pm^{1/N},
\end{equation}
where $\beta_\pm$ are roots of \eqref{e:80 N}.
Note that $0<\beta_-<\frac{N-1}{N}<\beta_+<1$ and hence
the solutions $\psi_\pm\in(0,\pi_N)$ of \eqref{e:B30} are uniquely (modulo $2\pi_N$) determined  and satisfy
$$0<\psi_+<\frac{\pi_N}{2}<\psi_-<\pi_N.$$
Observe that $F(\psi)$ is smooth, bounded and nonnegative for all $\psi\in\R$ and
strictly positive for $\psi\in(0,\psi_+)$.
Let $\psi(r)$ be the solution to the problem
\begin{equation}\label{e-psi-pN}
\psi'=\frac{{}(\psi)}{r\log r},\qquad\psi(R)=0,
\end{equation}
in $(R,+\infty)$,
for some $R>e
$.
Note also that $\psi_+(r)\equiv\psi_+$ is a stationary solution to  \eqref{e-psi-pN}.
So, $\psi(r)\le\psi_+$ for all $r>R$ by Lemma \ref{l:compar}.
Thus, by Lemma \ref{l:stab} we conclude that $\lim_{r\to\infty}\psi(r)=\psi_+$.

\begin{lemma}\label{l:B4pN} Let $\psi$ be the solution to \eqref{e-psi-pN}. Then
$\psi(r)$ admits a representation $\psi(r)=\psi_++\omega(r)$ where $\omega(r)<0$ in $[R,+\infty)$ and
\begin{eqnarray*}
\omega=c (\log r)^{\left(N(1-\beta_+)-1\right)(1+o(1))}
\quad\mbox{as}\quad r\rightarrow+\infty,\end{eqnarray*}
for some $c<0$.
\end{lemma}
The proof is the literary repetition of the arguments in the proof of Lemma \ref{l:B2}. Note only that $F(\psi_+)=0$,
$F_N'(\psi_+)=(1-\beta_+)N-1$.

Given the solution $\psi(r)$ to \eqref{e-psi-pN},
let $\rho(r)$ be the solution to the problem
\begin{equation}\label{e-rho-pN}
\frac{\rho'}{\rho}=\frac{G(\psi)}{r\log r}\quad\mbox{in}\quad(R,+\infty),\qquad\rho(R)=1.
\end{equation}
Observe that the right hand side of \eqref{e-rho-pN} is bounded and smooth
for all $(r,\psi)\in(R,\infty)\times\R$, so $\rho(r)$ exists for all $r>R$.

\begin{lemma}\label{p:CHeN} Let $\rho$ be the solution to \eqref{e-rho-pN}. Then
$\rho(r)=c(\log r)^{(\beta_+-1)(N-1)(1+o(1))}$ as ${r\rightarrow+\infty}$, for some $c>0$.
\end{lemma}
The proof is the literary repetition of the arguments in Lemma \ref{p:mu>0}. Notice only that
$G(\psi_+)=(1-\beta_+)(N-1)$.

\paragraph{(iv) Case $\mu=C_H$, $\epsilon \in(0,C_*)$, $p\neq N$.}
We consider in detail only the case $p>N$,
the case $p<N$  being similar.

The equations in system \eqref{e:BpsiV>0} can be written in the form
\begin{equation}\label{e:Bpsi_e}
\psi'=\frac{F_\eps(r,\psi)}{r},\qquad\frac{\rho'}{\rho}=\frac{G_\eps(r,\psi)}{r}\quad\mbox{in}\quad (R,+\infty),
\end{equation}
where we use the notation
$$
F_\eps(r,\psi):=
U(r)+W(r)S_p(\psi)S'_p(\psi)|S'_p(\psi)|^{p-2},\qquad
G_\eps(r,\psi):=W(r)S_p^{p}(\psi),$$
$$U(r):=\left(C_H+\frac{\epsilon}{\log^2r}\right)^{1/p},\qquad
W(r):=\frac{N-p}{p-1}-\frac{2\epsilon}{p\,(C_H\log^{2}r+\epsilon)\log r}.$$
Observe that $\frac{U(r)}{W(r)}\neq const$,
so the first equation in \eqref{e:BpsiV>0} has no stationary solutions.
For $\beta>0$, denote
\begin{eqnarray*}
A(r):=\gamma_\ast+\frac{\beta}{\log r}.
\end{eqnarray*}
Below we suppress the dependence on $r$ in $U(r)$ and $A(r)$ writing simply $U$ and $A$.
For $r> e$, let $\psi_\beta(r)$ be defined as the solution to the system
\begin{equation}\label{e:B38}
S_p(\psi)=\frac{U}{\left(|A|^p+\frac{U^p}{p-1}\right)^{1/p}}
\quad\mbox{and}\quad
S'_p(\psi)=\frac{A}{\left(|A|^p+\frac{U^p}{p-1}\right)^{1/p}},
\end{equation}
satisfying $0\leq \psi_\beta<\frac{\pi_p}{2}$.
From the definition of $U$ and $A$ one can see that $\lim_{r\to\infty}\psi_\beta(r)=(\pi/4)_p$,
where ${(\pi/4)_p}$ is defined by \eqref{e:B15}.

\begin{lemma}\label{l:B5} Let $\beta>0$. The following assertions are valid.
\begin{enumerate}
\item[$(i)$]If $\beta\in(\beta_-,\beta_+)$ then  there exists $R_\beta>e$ such that
 $\psi_\beta(r)$ is a positive super-solution to the equation
\begin{eqnarray}\label{e:psi_e'}
\psi'=\frac{F_\eps(r,\psi)}{r}\quad\mbox{in}\quad (R_\beta,+\infty).
\end{eqnarray}
\item[$(ii)$]If $\beta\not\in[\beta_-,\beta_+]$ then there exists $R_\beta>e$ such that
$\psi_\beta(r)$ is a  sub-solution
to  \eqref{e:psi_e'}.
\end{enumerate}
\end{lemma}
\begin{proof}

A routine calculation based on \eqref{e:B38} gives that
$$\psi_\beta'=\frac{1}{r}\,\left(F_\eps(r,\psi_\beta)+\frac{U}{A^p+\frac{U^p}{p-1}}\,\Theta(r)\right),$$
where
\begin{eqnarray*}
\Theta(r)&:=&
A^{p-2}\frac{\beta}{\log^2 r}-A^p-\frac{U^p}{p-1}-\frac{N-p}{p-1}A^{p-1}\\
&=&
\frac{\beta}{\log^2 r}\,\left|\gamma_\ast+\frac{\beta}{\log r}\right|^{p-2}
+\left(\frac{\gamma_\ast}{p-1}-\frac{\beta}{\log r}\right)\,\left|\gamma_\ast+\frac{\beta}{\log r}\right|^{p-1}
-\frac{\epsilon}{p-1}\,\frac{1}{\log^2 r}
-\frac{|\gamma_\ast|^p}{p-1}.
\end{eqnarray*}
For $r\to +\infty$ we obtain
$$\Theta(r)=-\frac{|\gamma_\ast|^{p-2}}{2p}\,
\frac{(\beta-\beta_+)(\beta-\beta_-)}{\log^2(r)}
+o\left(\frac{1}{\log^3 r}\right).
$$
Thus the assertion follows.
\end{proof}

Set $b:=\beta_+-\delta$, $B:=\beta_++\delta$, where $\delta>0$ is chosen such that
$\frac{1}{p}<b<\beta_+<B<\frac{2}{p}.$ By Lemma~\ref{l:B5} there exist $R_b$ and $R_B$ such that $\psi_b$
and $\psi_B$ are sub- and super-solution to \eqref{e:psi_e'}, respectively.
Set $R_\delta:=\max\{R_b,R_B\}$. It follows from \eqref{e:B38} that $\psi_B(R_\delta)<\psi_b(R_\delta)$.
Let $\psi_\ast(r)$ be the solution to the problem
\begin{eqnarray}\label{e-psi-ast}
\psi_\ast'=\frac{{F_\eps}(r,\psi_\ast)}{r}\quad\mbox{in}\quad(R_\delta,+\infty),\qquad\psi_\ast(R_\delta)=\psi_0,
\end{eqnarray}
where $\psi_0\in(\psi_B(R_\delta),\psi_b(R_\delta))$.
Observe that $F_\eps(r,\psi)$ is smooth and  bounded, so $\psi_\ast$ exists  for all $r>R_\delta$.
Moreover, by 
Lemma \ref{l:compar} we conclude that
\begin{eqnarray}\label{B-bound}
\psi_{B}(r)\leq\psi_\ast(r)\leq \psi_{b}(r),\qquad r\in[R_\delta,+\infty),
\end{eqnarray}
and, one can see that $F_\varepsilon(r,\psi_*(r))>0$ in $[R_\delta,+\infty)$.
Observe also that $\lim_{r\to\infty}\psi_\ast(r)=(\pi/4)_p$.

Let $\psi(r)$ be the solution to the problem
\begin{equation}\label{e-psi-eps}
\psi'=\frac{{F_\eps}(r,\psi)}{r},\qquad\psi(R_\delta)=0.
\end{equation}
Clearly, $\psi(r)$ exists  for all $r>R_\delta$. By Lemma \ref{l:compar} one has $0\leq\psi(r)\leq \psi_*(r)$.
Hence using the definitions of $F_\varepsilon$, $S_p$ and $S_p'$ one can see that
$F_\eps(r,\psi(r))$ is strictly positive
in $[R_\delta,+\infty)$.
Notice that
$$\lim_{r\to\infty}F_\eps(r,\psi)=F(\psi),$$
uniformly in $\psi$, where $F$ is defined by \eqref{e:FG-CH}.
Thus, by Lemma \ref{l:stab} we conclude that $\lim_{r\to\infty}\psi(r)=(\pi/4)_p$.

\begin{lemma}\label{l:B9} Let $\psi$ be the solution to \eqref{e-psi-eps} and $\psi_*$
be the solution to \eqref{e-psi-ast}. Then $\omega(r):=\psi(r)-\psi_*(r)<0$ in $[R_\delta,+\infty)$ and satisfies
the inequality
\begin{eqnarray*}
c_1\omega_+(r)\leq\omega(r)\leq c_2\omega_-(r),
\end{eqnarray*} for some  $c_1<0$, $c_2<0$, where
$\omega_+(r)=(\log r)^{-bp(1+o(1))}$, $\omega_-(r)=(\log r)^{-Bp(1+o(1))}$ as $r\to+\infty$.

\end{lemma}
\begin{proof}
Note that $\omega(R_\delta)=-\psi_0$ and $\omega(r)\to 0$ as $r\to +\infty$.
Fix $r> R_\delta$. Near $\psi_*(r)$ we have
\begin{eqnarray}\label{e:B44}
F_\eps(r,\psi)&=& F_\eps(\psi_*)+({F_\eps})_\psi'(r,\psi_*)(\psi-\psi_*)+
\frac{1}{2}({F_\eps})_\psi''(r,\psi_*)(\psi-\psi_*)^2+\Theta(r,\psi-\psi_*)
\end{eqnarray}
where $\Theta(r,\psi-\psi_*)=o((\psi-\psi_*)^2)$ as $\psi\rightarrow \psi_*$.
A direct computation gives
\begin{eqnarray}\label{e:B42}
(F_\eps(r,\psi))'_\psi&=&
W(r)\,(|S_p'(\psi)|^p-|S_p(\psi)|^p);\\
\label{e:B43}
(F_\eps(r,\psi))''_\psi&=&
-\frac{p^2}{p-1}\,W(r)|S_p(\psi)|^{p-1}S_p'(\psi).
\end{eqnarray}
Since $\psi=\psi_*+\omega$, and $\psi_*$ solves the same equation, from \eqref{e:B44} we obtain
\begin{equation}\label{e:B-2}\omega'=({ F_*})_\psi'(\psi_*)\frac{\omega}{r}+\frac{({ F_*})_\psi''(\psi_*)}{2}\,\frac{\omega^2}{r}
+ \frac{\Psi(\omega)}{r}.\end{equation}
Using \eqref{B-bound}, \eqref{e:B38}, \eqref{e:B42} and \eqref{e:B43} we conclude that
$$-\frac{Bp}{\log r}+o\left(\frac{1}{\log^2 r}\right)
\leq (F_*)_\psi'(\psi_*(r))\leq
-\frac{bp}{\log r}+o\left(\frac{1}{\log^2 r}\right),
$$
$$
\frac{bp^2}{p-1}\,\frac{1}{\log r}+o\left(\frac{1}{\log^2 r}\right)
\leq (F_*)_\psi''(\psi_*(r))-\frac{(p-N)p}{p-1}\leq
\frac{Bp^2}{p-1}\,\frac{1}{\log r}+o\left(\frac{1}{\log^2 r}\right),
$$
as $r\rightarrow +\infty$.
We substitute the above estimates   into \eqref{e:B-2}. Then
\begin{equation}\label{e:B45}
\omega'\leq
-Bp\frac{\omega}{r\log r}
+\frac{(p-N)p}{2(p-1)}\frac{\omega^2}{r}
+\frac{Bp^2}{2(p-1)}\frac{\omega^2}{r\log r}
+O\left(\frac{\omega+\omega^2}{r\log^2 r}+\frac{\omega^3}{r}\right),
\end{equation}
\begin{equation}\label{e:B45-2}
\omega'\geq
-bp\frac{\omega}{r\log r}
+\frac{(p-N)p}{2(p-1)}\frac{\omega^2}{r}
+\frac{bp^2}{2(p-1)}\frac{\omega^2}{r\log r}
+O\left(\frac{\omega+\omega^2}{r\log^2 r}+\frac{\omega^3}{r}\right),
\end{equation}
as $r\rightarrow +\infty$.
From \eqref{e:B45-2} we infer that
$$\frac{\omega'}{\omega}\leq
-\frac{bp}{r\log r}
+\frac{bp^2}{p-1}\,\frac{\omega}{r\log r}
+O\left(\frac{1+\omega}{r\log^2 r}+\frac{\omega^2}{r}\right),$$
and, hence,
$$\log\frac{\omega(r)}{\omega(R)}\leq c\log (\log r)^{-bp(1+o(1))},$$
or, equivalently,
$$\omega(r)\geq
-c(\log r)^{-bp(1+o(1))},$$
as $ r\rightarrow+\infty$.
Therefore by \eqref{e:B45} we infer that $\omega$ satisfies
$$
\frac{\omega'}{\omega}\geq
-\frac{Bp}{r\log r}
-\frac{c}{r(\log r)^{bp(1+o(1))}}
+\frac{Bp^2}{2(p-1)}\,\frac{\omega}{r\log r}
+O\left(\frac{1+\omega}{r\log^2 r}+\frac{\omega^2}{r}\right),
$$
and, hence,
$$\log\frac{\omega(r)}{\omega(R)}\geq\log (\log r)^{-Bp(1+o(1))}
-c(\log r)^{-bp(1+o(1))+1},$$
as $ r\rightarrow+\infty$.
Since $b>\frac{1}{p}$,  one has
$$\omega(r)\leq
-c(\log r)^{-Bp(1+o(1))},$$
as $r\rightarrow +\infty$.
The assertion follows.
\end{proof}

Given the solution $\psi(r)$ to the problem \eqref{e-psi-eps},
let $\rho(r)$ be the solution to the problem
\begin{equation}\label{e-rho-eps}
\frac{\rho'}{\rho}=\frac{G_\eps(r,\psi)}{r}\quad\mbox{in}\quad(R_\delta,+\infty),\qquad\rho(R_\delta)=1.
\end{equation}
Clearly the right hand side of \eqref{e-rho-eps} is bounded and smooth
for all $(r,\psi)\in(R_\beta,\infty)\times\R$, so $\rho(r)$ exists  for all $r>R_\beta$.

\begin{lemma}\label{p:Che}
 Let $\rho$ be the solution to \eqref{e-rho-eps}.
Then $\rho$ satisfies the estimate
\begin{eqnarray*}
c_1 \rho_-(r)\leq\rho(r)
\leq c_2 \rho_+(r)
\end{eqnarray*}
for some $c_1>0$, $c_2>0$, where $\rho_-(r)=r^{-\gamma_\ast}(\log r)^{b(p-1)(1+o(1))}$,
$\rho_+(r)=r^{-\gamma_\ast}(\log r)^{B(p-1)(1+o(1))}$  as $r\to+\infty$.
\end{lemma}

\begin{proof}
Near $\psi_*(r)$ we have
$G_\eps(r,\psi)=G_\eps(r,\psi_*)+\Xi(r,\psi-\psi_*)$,
where $\Xi(r,\psi-\psi_*)=O(\psi-\psi_*)$ as $r\rightarrow +\infty$.
Using \eqref{B-bound} 
 we obtain
$$
W(r)S_p^{p}(\psi_{b})\leq G_\eps(r,\psi_*)\leq W(r)S_p^{p}(\psi_{B}).
$$
By a simple computation we conclude that
$$\frac{b(p-1)}{\log r}+\frac{o(1)}{\log^2r}
\leq G_\eps(r,\psi_*)+\gamma_\ast\leq
\frac{B(p-1)}{\log r}+\frac{o(1)}{\log^2r},$$
as $r\rightarrow+\infty$.
Therefore $\rho$ satisfies
$$\frac{b(p-1)}{r\log r}
+O\left(\frac{1}{r\log^2r}+\frac{\omega}{r}\right)
\leq\frac{\rho'}{\rho}+\frac{\gamma_\ast}{r}\le
\frac{B(p-1)}{r\log r}
+O\left(\frac{1}{r\log^2r}+\frac{\omega}{r}\right),
$$
as $r\rightarrow+\infty$.
Thus we have
\begin{eqnarray*}
&&\log\left(\frac{\log r}{\log R}\right)^{b(p-1)}
+\frac{O(1)}{\log r}+c(\log r)^{-bp(1+o(1))+1}
\\&&\leq
\log\frac{\rho(r)}{\rho(R)}-\log\left(\frac{r}{R}\right)^{-\gamma_\ast}
\leq\log\left(\frac{\log r}{\log R}\right)^{B(p-1)}
+\frac{O(1)}{\log r}+c(\log r)^{-Bp(1+o(1))+1},
\end{eqnarray*}
as $r\rightarrow+\infty$.
So, the assertion follows from $\frac{1}{p}<b<B$.
\end{proof}

\begin{remark}
The case $\mu=C_H$, $\epsilon\in(0,C_\ast)$ and $p<N$ is similar,
the only difference being 
$\gamma_*<0$ and hence $\lim_{r\to+\infty}\psi_\ast(r)=(3\pi/4)_p$.
\end{remark}

\end{small}

\begin{small}
\paragraph{Acknowledgments.}
The authors are grateful to Wolfgang Reichel for his suggestion to use the generalized
Pr\"ufer transformation for obtaining the asymptotic of large sub-solutions in Appendix C,
and to Friedemann Brock for his valuable comments.
The financial support of the Royal Society is gratefully acknowledged.
\end{small}

\begin{small}

\end{small}

\begin{thebibliography}{99}

\bibitem{Peral-1}
{\sc B.\,Abdellaoui and\ I.\,Peral}, {\it Existence and nonexistence
results for quasilinear elliptic equations involving the
$p$-Laplacian with a critical potential}.
Ann. Mat. Pura Appl. (4) {\bf 182} (2003), 247--270.

\bibitem{Peral-2}
{\sc B.\,Abdellaoui and\ I.\,Peral},
{\it On quasilinear elliptic equations related to some Caffarelli-Kohn-Nirenberg inequalities.}
Commun. Pure Appl. Anal. {\bf 2} (2003), 539--566.



\bibitem{AdimChandRamas}
{\sc Adimurthi, N.\,Chaudhuri and M.\,Ramaswamy}, {\it An improved
Hardy-Sobolev inequality and its application}. {Proc. Amer. Math.
Soc.} {\bf 130} (2002), 489--505.

\bibitem{Agmon}
{\sc S.\,Agmon}, {\it On positivity and decay of solutions of
second order elliptic equations on Riemannian manifolds}, in {\it
Methods of functional analysis and theory of elliptic equations
(Naples, 1982)}, 19--52, Liguori, Naples, 1983.

\bibitem{Agmon-2}
{\sc S.\,Agmon}, {\it Bounds on exponential decay of eignfunctions
of Schr\"odinger operators},  Schr\"odinger operators (Como, 1984),
1--38, Lecture Notes in Math., 1159, Springer, Berlin, 1985.



\bibitem{AllegrettoHuang}
{\sc W.\,Allegretto and Y.\,H.\,Huang},  {\it A Picone's identity for the
$p$-Laplacian and applications}. {Nonlinear Analysis TMA}
{\bf 32}(1998), 819--830.

\bibitem{Barbatis_Fillipas_Tertikas}
{\sc G.\,Barbatis, S.\,Filippas and A.\,Tertikas}, {\it A unified
approach to improved $L\sp p$ Hardy inequalities with best
constants}. {Trans. Amer. Math. Soc}. {\bf 356} (2004), 2169--2196.

\bibitem{B-Veron}
{\sc M.-F.\,Bidaut-V\'eron}, {\it Local and global behavior of
solutions of quasilinear equations of Emden-Fowler type}. Arch.
Rational Mech. Anal. {\bf 107} (1989), 293--324.

\bibitem{B-VP}
{\sc M.-F.\,Bidaut-V\'eron and S.\,Pohozaev},
{\it Nonexistence results and estimates for some nonlinear elliptic problems}.
J. Anal. Math. {\bf 84} (2001), 1--49.

\bibitem{Binding}
{\sc P.\,Binding and P.\,Drabek},
{\it Sturm-Liouville theory for the $p$-Laplacian.}
Studia Sci. Math. Hungar. {\bf 40} (2003), 375--396.


\bibitem{Boggio}
{\sc T.\,Boggio}, {\it Sull'equazione del moto vibrattorio delle
membrane elastiche}. {Accad. Lincei, Sci. Fis.}, ser. 5, {\bf 15} (1907), no. 6, 386--393.

\bibitem{Brezis-Cabre}
{\sc H.\,Brezis and X.\,Cabr\'e}, {\it Some simple nonlinear PDE's
without solutions}. {Boll. Unione Mat. Ital. Sez. B} (8) {\bf 1}
(1998), 223--262.

\bibitem{Brezis-Tesei}
{\sc H. Brezis, L. Dupaigne and A. Tesei},
{\it On a semilinear elliptic equation with inverse-square potential}.
Selecta Math. (N.S.) {\bf 11} (2005), 1--7.

\bibitem{Brown-Reichel}
{\sc B.\,M.\,Brown and W.\,Reichel},
{\it Eigenvalues of the radially symmetric $p$-Laplacian in $\mathbb R\sp n$}.
J. London Math. Soc. (2) {\bf 69} (2004), 657--675.

\bibitem{Cod}
{\sc E.\,A.\,Coddington and N.\,Levinson},
{\it Theory of ordinary differential equations.}
McGraw-Hill Book Company, Inc., New York-Toronto-London, 1955. xii+429 pp.

\bibitem{DelPino}
{\sc M.\,Del Pino, M.\,Elgueta and R.\,Man\'asevich},
{\it A homotopic deformation along $p$ of a Leray-Schauder degree result and existence
for $(|u'|^{p-2}u')'+f(t,u)=0$, $u(0)=u(T)>0$, $p>1$}.
{J. Differential Equations} {\bf 80} (1989), 1--13.


\bibitem{Diaz-Saa}
{\sc J.\,I.\,D\'iaz and J.\,E.\,Saa},
{\it Existence et unicit\'e de solution positive pour certaines \'equations elliptiques quasilin\'eaires},
C.R. Acad. Sci. Paris {\bf 305} (1987) 521--524.

\bibitem{Drabek}
{\sc P. Dr\'abek and J. Hern\'andez},
{\it Existence and uniqueness of positive solutions for some quasilinear elliptic problems.}
Nonlinear Anal. {\bf 44} (2001), 189--204.

\bibitem{Dupaigne}
{\sc L.\,Dupaigne}, {\it Semilinear elliptic PDE's with a singular
potential}, {J. Anal. Math.} {\bf 86} (2002), 359--398.

\bibitem{Barbatis}
{\sc  S.\,Filippas and A.\,Tertikas},
{\it Optimizing improved Hardy inequalities}.
{J. Funct. Anal.} {\bf 192} (2002), 186--233.

\bibitem{Fleck-1}
{\sc J.\,Fleckinger, J.\,Hern\'andez, P.\,Tak\'a\v c, and F.\,de\,Th\'elin},
{\it Uniqueness and positivity for solutions of equations with the $p$-Laplacian},
Proceedings of the Conference on Reaction-Diffusion Equations, 1995, Trieste, Italy.
Lecture Notes in Pure and Applied Math. Vol. 194 (1998), 141--155.
Marcel Dekker, New York and Basel.

\bibitem{Harrell}
{\sc J.\,Fleckinger, E.\,M.\,Harrell\,II and  F.\,de\,Th\'elin},
{\it Boundary behavior and  $L^q$ estimates for solutions of equations containing the $p$-Laplacian.}
Electron. J. Differential Equations  1999, No. 38, 1-19 p. (electronic).


\bibitem{Melian-Sabina de Lis}
{\sc J.\,Garc\'ia-Meli\'an and J.\,Sabina de Lis}, {\it Maximum and
Comparison Principles for Operators Involving the $p$-Laplacian}.
{ J. Math. Anal. Appl.} {\bf 218} (1998), 49--65.

\bibitem{GazzolaGrunauMitidieri}
{\sc F.\,Gazzola, H-C.\,Grunau and E.\,Mitidieri}, {\it Hardy
inequalities with optimal constants and remainder terms}. Trans.
Amer. Math. Soc. {\bf 356} (2004), 2149--2168.

\bibitem{Gidas-Spruck}
{\sc B.\,Gidas and J.\,Spruck}, {\it Global and local behavior of
positive solutions of nonlinear elliptic equations}, {Comm. Pure
Appl. Math.} {\bf 34} (1981), 525--598.



\bibitem{Kondratiev-Landis}
{\sc V.\,Kondratiev and E. Landis},
{\it Qualitative properties of the solutions of a second-order nonlinear equation.} (Russian)
Mat. Sb. (N.S.) {\bf 135}(177) (1988), 346--360, 415.

\bibitem{KLS}
{\sc V.\,Kondratiev, V.\,Liskevich and Z.\,Sobol},
{\it Second--order semilinear elliptic inequalities in exterior domains}.
J. Differential Equations {\bf 187} (2003), 429--455.

\bibitem{KLS1}
{\sc V.\,Kondratiev, V.\,Liskevich and Z.\,Sobol}, {\it Positive
super-solutions to semilinear second--order non-divergence type
elliptic equations in exterior domains}. Preprint, 2005.


\bibitem{KLM}
{\sc V.\,Kondratiev, V.\,Liskevich and V.\,Moroz}, {\it Positive
solutions to superlinear second--order divergence type elliptic
equations in cone--like domains}.
Ann. Inst. H. Poincar\'e Anal. Non Lin\'eaire 22 (2005), 25-43.

\bibitem{KLMS}
{\sc V.\,Kondratiev, V.\,Liskevich, V.\,Moroz and Z.\,Sobol}, {\sl A
critical phenomenon for sublinear elliptic equations in cone--like
domains}.
Bull. Lond. Maths. Soc. {\bf 37} (2005), 585-591.


\bibitem{Lindquist-Prufer}
{\sc P.\,Lindquist}, {\it Some remarkable sine and cosine
functions}. Ricerche Mat. {\bf 44} (1995), 269--290.

\bibitem{LLM}
{\sc V.\,Liskevich, S.\,Lyakhova and V.\,Moroz},
{\it Positive solutions to singular semilinear elliptic equations
with critical potential on cone-like domains}.
Adv. Differential Equations. {\bf 11} (2006), 361-398.

\bibitem{LiskevichSkrypnikSkrypnik}
{\sc V.\,Liskevich, I.I.\,Skrypnik and I.V.\,Skrypnik} {\it
Positive super-solutions to general nonlinear elliptic equations
in exterior domains}. Manuscripta Math. {\bf 115} (2004), no. 4, 521--538.


\bibitem{ShafrirMarcus}
{\sc M.\,Marcus and I.\,Shafrir}, {\it An eigenvalue problem related to
Hardy's $L\sp p$ inequality}.  {Ann. Scuola Norm. Sup. Pisa Cl.
Sci}. (4) {\bf 29} (2000), 581--604.

\bibitem{Pohozhaev}
{\sc E.\,Mitidieri and S.\,I.\,Poho\v zaev}, {\it  A priori
estimates and the absence of solutions of nonlinear partial
differential equations and inequalities (Russian)}. Tr. Mat. Inst.
Steklova {\bf 234} (2001), 1--384.


\bibitem{Pinchover-Tintarev}
{\sc Y.\,Pinchover and  K.\,Tintarev},
{\it Ground state alternative for $p$-Laplacian with potential term},
Preprint, 2005.

\bibitem{Pohozaev-Tesei}
{\sc S. I. Pohozaev\ and\ A. Tesei},
{\it  Nonexistence of local solutions to semilinear partial differential inequalities},
Ann. Inst. H. Poincar\'e Anal. Non Lin\'eaire {\bf 21} (2004), 487--502.

\bibitem{ShafrirPoljakovsky}
{\sc A.\,Poliakovsky and I.\,Shafrir}, {\it A comparison principle for
the $p$-Laplacian. Elliptic and parabolic problems} (Rolduc/Gaeta,
2001), 243--252, {World Sci. Publishing, River Edge}, NJ, 2002.

\bibitem{ReichelWalter}
{\sc W.\,Reichel and W.\,Walter}, {\it Sturm-Lioville Type Problems
for the $p$-Laplacian under Asymptotic Non-resonance Conditions}.
J. Differential Equations  {\bf 156}, (1999), 50--70.

\bibitem{Sandeep}
{\sc K.\,Sandeep}, {\it On the first eignfunction of a perturbed
Hardy-Sobolev operator}.   NoDEA Nonlinear Differential Equations Appl. {\bf 10}
(2003), 223--253.

\bibitem{SandeepSreenadh}
{\sc K.\,Sandeep and K.\,Sreenadh}, {\it Asymptotic behavior of the
first eignfanction of a Hardy-Sobolev operator}.  Nonlinear
Anal. {\bf 54} (2003), 545--563.

\bibitem{Shafrir}
{\sc I.\,Shafrir}, {\it Asymptotic behavior of minimizing sequences
for Hardy's inequality}. Commun. Contemp. Math., {\bf 2} (2002),
151--189.

\bibitem{Serrin64}
{\sc J.\,Serrin}, {\it Local behavior of solutions of
quasi-linear elliptic equations}, Acta Math. {\bf 111} (1964),
247--302.

\bibitem{Serrin-Zou}
{\sc J.\,Serrin and H.\,Zou}, {\it Cauchy--Liouville and universal
boundedness theorems for quasilinear elliptic equations and
inequalities}, Acta Math. {\bf 189} (2002), 79--142.

\bibitem{Smets}
{\sc D.\,Smets and A.Tesei},
{\it On a class of singular elliptic problems with first order terms},
{Adv. Differential Equations} {\bf 8} (2003), 257--278.

\bibitem{Takac}
{\sc P.\,Tak\'a\v c, L.\,Tello and\ M.\,Ulm}, {\it Variational
problems with a $p$-homogeneous energy}. Positivity {\bf 6}
(2002), 75--94.

\bibitem{Terracini}
{\sc S.\,Terracini}, {\it On positive entire solutions to a class of
equations with a singular coefficient and critical exponent}. Adv.
Differential Equations {\bf 1} (1996), 241--264.


\bibitem{Veron}
{\sc L. V\'eron},
 Singularities of solutions of second order quasilinear equations.
Pitman Research Notes in Mathematics Series, 353.
Longman, Harlow, 1996. viii+377 pp.



\end{thebibliography}
\end{document}